\documentclass[times,sort&compress,3p]{elsarticle}
\journal{Journal of Multivariate Analysis}
\usepackage[labelfont=bf]{caption}

\usepackage{amsmath,amsfonts,amssymb,amsthm,booktabs,color,epsfig,graphicx,url}

\RequirePackage[colorlinks,citecolor=blue,urlcolor=blue]{hyperref}
\RequirePackage{graphicx,subfigure}
\usepackage{bbm,bm}
\usepackage{enumerate}
\usepackage{caption}
\usepackage{xcolor}
\usepackage{threeparttable}
\usepackage{makecell}

\theoremstyle{plain}
\newtheorem{thm}{Theorem}[section]

\newtheorem{lemma}{Lemma}[section]
\newtheorem{cor}{Corollary}[section]
\theoremstyle{definition}

\newtheorem{remark}{Remark}[section]

\newcommand{\bbV}{{\bf V}}
\newcommand{\bbU}{{\bf U}}
\newcommand{\bbD}{{\bf D}}
\newcommand{\bbA}{{\bf A}}

\newcommand{\bbX}{{\bf X}}
\newcommand{\bbP}{{\bf P}}

\newcommand{\bbS}{{\bf S}}
\newcommand{\bbB}{{\bf B}}

\newcommand{\bSi}{\pmb \Sigma}

\newcommand{\bgl}{{\bf \lambda}}

\newcommand{\bGma}{{\bf \Gamma}}

\newcommand{\bgL}{{\bf \Lambda}}

\newcommand{\bUps}{{\bf \Upsilon}}

\newcommand{\bbI}{{\bf I}}

\newcommand{\bbT}{{\bf T}}

\newcommand{\bqn}{\begin{eqnarray*}}
	\newcommand{\eqn}{\end{eqnarray*}}

\newcommand{\bqa}{\begin{eqnarray}}
	\newcommand{\eqa}{\end{eqnarray}}

\newcommand{\al}{\alpha}

\newcommand{\CYRS}{\CYRS}

\begin{document}

\begin{frontmatter}

\title{Asymptotic distributions of four linear hypotheses test statistics under generalized spiked model}

\author[1]{Zhijun Liu}
\author[2]{Jiang Hu\corref{mycorrespondingauthor}}
\author[2]{Zhidong Bai}
\author[4]{Zhihui Lv}

\address[1]{College of Sciences, Northeastern University, 
	China}
\address[2]{KLASMOE and School of Mathematics and Statistics, Northeast Normal University, 
China	}
\address[4]{School of Mathematics and Statistics,  Guangdong University of Foreign Studies, 
	China}
\cortext[mycorrespondingauthor]{Corresponding author. Email address: \url{huj156@nenu.edu.cn}}

\begin{abstract}
In this paper, we establish the Central Limit Theorem (CLT) for linear spectral statistics (LSSs) of large-dimensional generalized spiked sample covariance matrices, where the spiked eigenvalues may be either bounded or diverge to infinity. Building upon this theorem, we derive the asymptotic distributions of linear hypothesis test statistics under the generalized spiked model, including Wilks’ likelihood ratio test statistic $U$, the Lawley-Hotelling trace test statistic $W$, and the Bartlett-Nanda-Pillai trace test statistic $V$. Due to the complexity of the test functions, explicit solutions for the contour integrals in our calculations are generally intractable. To address this, we employ Taylor series expansions to approximate the theoretical results in the asymptotic regime. We also derive asymptotic power functions for three test criteria above, and make comparisons with Roy’s largest root test under specific scenarios. Finally, numerical simulations are conducted to validate the accuracy of our asymptotic approximations.
\end{abstract}


\end{frontmatter}

\section{Introduction}
\label{result UWV}
Linear hypothesis testing plays an important role in the analysis of multivariate data. Four criteria used to test linear hypotheses are: Wilks' likelihood ratio criterion, Lawley-Hotelling trace criterion,  Bartlett-Nanda-Pillai trace criterion, and Roy's largest root criterion. The corresponding test statistics are defined as:
\begin{itemize}
	\item Wilks' likelihood ratio $U=\sum_{i=1}^{p}\log\left(1+\lambda_{i} \right)  $
	\item Lawley-Hotelling trace $W=\sum_{i=1}^{p}\lambda_{i} $
	\item Bartlett-Nanda-Pillai trace $ V=\sum_{i=1}^{p}\frac{\lambda_{i}}{1+\lambda_{i}} $
	\item Roy's largest root $R=\lambda_{1} $
\end{itemize}	
where $ \lambda_{i}, i=1,\dots,p $ are the eigenvalues of an $ F $ matrix, which is the product of a sample covariance matrix from the independent variable array $ (x_{ij})_{p\times n_{1}} $ and the inverse of
another covariance matrix from the independent variable array $ (y_{ij})_{p\times n_{2}} $. Based on the differences of the four test functions, we divide the four statistics into two categories. The first category includes statistics which are extreme eigenvalues of a matrix, such as the largest eigenvalue or the smallest eigenvalue. For example, statistic $R$ belongs to this category. In the second category, statistics can be expressed as a linear combination of the function of all the eigenvalues, such as $U, W, V$, and they are also called linear spectral statistics (LSSs). 

In this work, we consider the general sample covariance matrix  $ \bbB_n=\frac{1}{n}\bbT_p\bbX_n\bbX_n^{\ast}\bbT_p^{\ast} $, where $ \bbX_n $ is a $ p\times n $ matrix with independent and identically distributed (i.i.d.) standardized entries $ \left\lbrace x_{ij}\right\rbrace _{1\leq i\leq p, 1\leq j \leq n} $, $ \bbT_p $ is a $p \times p$ deterministic matrix, $\bbT_p\bbX_n$ is considered a random sample from the population  covariance matrix $\bbT_p\bbT_p^{\ast}=\bSi$, and $^*$ represents
the complex conjugate transpose. In the sequel, we simply write $\bbB\equiv\bbB_n$, $\bbT\equiv\bbT_p$ and $\bbX\equiv\bbX_n$ when there is no confusion. We denote $ \lambda_1\geq \lambda_2\geq \dots\geq \lambda_p$ as eigenvalues of $\bbB$, and denote $\alpha_1\geq \al_2\geq \dots\geq \al_p$ as the eigenvalues of $\bSi$. For a known
test function $f$, we call $\sum_{j=1}^{p}f(\lambda_j)$ an LSS of $\bbB$. For example, for three LSSs $U,W,V$, test functions are $f_U=\log(1+x)$, $f_W=x$, $f_V=x/(1+x)$, respectively.

For the aforementioned first category, the primary interest lies in the asymptotic behavior of a few largest eigenvalues and their eigenvectors.  
According to the seminal work of \cite{Baik05}, we have that the largest eigenvalue of sample covariance $\bbB$ undergo a phase transition: define $c=p/n$, when  $\al_1<1+\sqrt{c}$, $\lambda_1$  converges to the right end of Marchenko-Pastur law (MP law); when $\al_1>1+\sqrt{c}$, eigenvalue $\lambda_1$ jumps out of the support of MP law; Moreover, phase transition also happens on the level of the second-order fluctuation. Specifically, when $ \alpha_1-(1+\sqrt{c})\ll n^{-\frac{1}{3}}$, which is also called subcritical regime. Under subcritical regime, $\lambda_1$ admits Tracy-Widom distribution; when $ \alpha_1-(1+\sqrt{c})\gg n^{-\frac{1}{3}}$, which is also called supercritical regime, then $\lambda_1$ has an asymptotic Gaussian distribution; when  $ \alpha_1-(1+\sqrt{c})\sim n^{-\frac{1}{3}}$, which is also called critical regime, then $\lambda_1$ has an asymptotic distribution between Tracy-Widom distribution and Gaussian distribution.

For the second  category, many efforts have been put into the  properties of LSSs under the high-dimensional case. As a benchmark, \cite{10.1214/aop/1078415845}  established the central limit theorem (CLT) for the LSSs of a high-dimensional $ \bbB $ under Gaussian-like moments condition by employing random matrix theory (RMT). Here the term `Gaussian-like moments' refers to the population second-order and fourth-order moments are the same as those of real or complex standard normal distribution. Following the work of \cite{10.1214/aop/1078415845}, many extensions have been developed under many different settings. For example, {\cite{pan2008central} relaxed the Gaussian-like moments condition of $ x_{ij} $, which  added a structural condition on $ \bbT $.  } \cite{LP09}, \cite{zheng2012central} and \cite{BaiH15C}  extended the BST to multivariate Wigner matrices, $F$ matrices and Beta matrices, respectively. \cite{YangP15I}, \cite{ZhengC19T} and \cite{BaoH22S} focused on the LSS for CLT of canonical correlation matrices, correlation matrices and block correlation matrices. \cite{GaoH17H} and \cite{LiW21C} studied the CLTs for the LSSs of high-dimensional Spearman and Kendall's rank correlation  matrices, respectively. \cite{pan2014comparison} presented the CLT for the LSS of noncentered sample covariance matrices, and \cite{10.1214/14-AOS1292} studied the case of an unbiased sample covariance matrix when the population mean is unknown. 
Under the ultra-high dimensional setting, \cite{chen2015clt} focused on the ultrahigh dimensional case in which the dimension $p$ is much larger than the sample size $n$. Compared with \cite{chen2015clt}, \cite{Qiu23} studied a more general setting, they considered a re-normalized sample covariance matrix and derived the asymptotic normality for spectral statistics of  the re-normalized sample covariance matrix when $p/n\rightarrow\infty$. Without attempting to be comprehensive, we also refer readers to other extensions \citep{BaiM07A, BaiL19C, BannaN20C, najim2016gaussian, baik2018ferromagnetic, hu2019high, li2020asymptotic, zhang2022asymptotic}. 

However, almost all the literature about asymptotic distributions of LSSs have traditionally assumed that the population covariance matrices are bounded in $n$, and this assumption cannot be satisfied in certain fields, such as signal detection or factor model, see examples in \cite{Liu23}. Recently, under the unbounded population setting, \cite{yin2021spectral} investigated the asymptotic distribution of LSS for sample covariance matrices when test function $f=x,x^2$.	Some other investigations about the unbounded population setting can be found in \cite{Liu23}; \cite{zhao23}; \cite{yin24}; \cite{Li25}.
\textcolor{black}{In this paper, we focus on a generalized spiked covariance model, which is defined as 
\begin{align}\label{ds}
	\bSi=\mathbf{V}\left(\begin{array}{cc}
		\bbD_{1} & 0 \\
		0 & \bbD_{2}
	\end{array}\right) \mathbf{V}^{\ast},
\end{align}
where $\mathbf{V}$ is a unitary matrix, $\bbD_{1}$ is a diagonal matrix with its elements are the spiked eigenvalues of $\bSi$, and they can be bounded or tend to infinity, and $\bbD_{2}$ is the diagonal matrix of the bulk eigenvalues.  Compared with \cite{Liu23}, model (\ref{ds}) is more general since spiked eigenvalues in (\ref{ds}) can be bounded spikes. Consequently, this model has a wider scope of application. To provide explicit formulas for the asymptotic means and variances for $U,W,V$ (Theorems \ref{U}--\ref{V}), we assume $\bbD_{2}$ is an identity matrix in Section \ref{section4}. In Section \ref{proof of UWV}, we obtain a CLT (Theorem \ref{thm1}) when $\bSi$ is a generalized spiked covariance model (\ref{ds}).}
Actually, model (\ref{ds}) is attributed to the famous spiked model proposed by \cite{10.1214/aos/1009210544}, in which a few large eigenvalues of the population covariance matrix are assumed to be well separated from the remaining eigenvalues. The spiked model has served as the foundation for a rich theory of principal component analysis through the performance of extreme eigenvalues, and signicant progress has been made on this topic in the recent few years, as discussed in \cite{BAIK20061382, 10.2307/24307692, 10.1214/07-AIHP118, Nadler08F,JungM09P,  BAI2012167, OnatskiM14S, BloemendalK16P, WangY17E, DonohoG18O, JohnstoneP18P, JohnstoneO20T, cai2020limiting, JiangB21G}.

In some sense, since statistics in the aforementioned first category are based on part of eigenvalues, therefore they are also called local statistics, whereas  statistics in the second category are based on all the eigenvalues,  then they are also called global statistics.
Comparisons between local statistics and global statistics have consistently attracted significant attention from researchers. To name a few, \cite{Olson74} concerned with the fixed-effects model of multivariate analysis of variance and compared $ U, W, V, R$ with other two test criteria by using Monte Carlo methods.
\cite{Dobriban16} concluded that tests based on top eigenvalue alone have small power to detect weak signals in high dimensions, therefore, to detect weak signals, an optimal inference should be based on all eigenvalues.
\cite{Ding23} analyzed superior power for global statistics and local statistics under general local alternatives when dimension $p$ is much larger than sample size $n$.
Recently, \cite{Liu23} compared the corrected likelihood ratio test and corrected Nagao's trace test with Roy's largest root test under the spiked model when the number of spikes is not always equal to $1$.

In this work, we obtain a generalized CLT for the LSSs of sample covariance matrices under population (\ref{ds}),
and the established CLT is employed to study the asymptotic behaviors of test statistics $U, W, V$  under the hypothesis 
\begin{align} \label{hytest1}
	H_{0}:\bSi=\bbI_{p}\quad \text{v.s.} \quad
	H_{1}:\bSi=\mathbf{V}\left(\begin{array}{cc}
		\bbD_{1} & 0 \\
		0 & \bbI_{p-M}
	\end{array}\right) \mathbf{V}^{\ast}.
\end{align} 
Because of the complexity of the test function, explicit solutions for the contour integrals in our calculations are generally intractable. To address this, we employ Taylor series expansions to approximate the theoretical results in the asymptotic regime.
Numerical simulations further confirm that our asymptotic results provide a highly accurate approximation. 
We also derive the asymptotic powers of four tests to detect hypothesis (\ref{hytest1}) and make comparisons between them under certain scenarios. We now describe the main contributions of the present paper as follows.
\textcolor{black}{
\begin{itemize}
	\item 
	First, compared with the traditional computations of the asymptotic mean and variance of LSSs, our approach introduces methodological innovations. When test functions of LSSs are $x,\log(x)$ or their linear combinations, one can use the residue theorem directly to calculate the asymptotic mean and variance, such as \cite{li2020asymptotic,WangY17E}. However, when test functions are complex, such as $f_U,f_V$, it is difficult to use the residue theorem directly. In this paper, we employ Taylor expansions to test functions $f_U,f_V$ and finally we use Taylor series expansions to approximate the theoretical results in the asymptotic regime.	In Section \ref{section5}, we provide some simulations to check the accuracy of the theoretical results.
	\item For fixed dimension, \cite{Johnstone17} developed an accurate and tractable
	asymptotic distribution of $R$ under a rank-one alternative, which is a combination of central and noncentral $\chi^2$ and $F$ variates with the restriction of divergent parameter. In high-dimensional case, \cite{hou19} applied the spiked model theory to develop a new method to obtain the asymptotic distribution of $R$ under a rank-finite alternative. Although \cite{Johnstone17} and \cite{hou19} also mentioned the test statistics $U, W, V$, they predominantly concentrated on the $R$. In this work, in addition to the $R$, we also take into account three other linear hypothesis test statistics, $ U, W,$ and $V$, and  we also make comparisons between them. Compared with classical works that study the $U, W, V,$ and $R$ test statistics such as \cite{Olson74}, our results are obtained under high-dimensional settings and do not require the normality assumptions.
	\item Compared with \cite{Liu23}, there are some differences in model setting and final asymptotic results. For the model setting,  a more reasonable approach is to place all spiked eigenvalues in the matrix $\bbD_1$. Moreover, we allow spiked eigenvalues in $\bbD_1$ to be bounded or diverge to infinity, but not all the spiked eigenvalues are diverging spikes. In \cite{Liu23}, they derive the asymptotic distributions of two common test statistics $\mathrm{tr}\bbB-\log\left| \bbB\right|-p $ and $ \mathrm{tr}(\bbB-\bbI_p)^2$. In this work, we consider four different statistics $U, W, V, R$ and perform thorough comparisons of four tests based on $U, W, V, R$. 
\end{itemize}
}

The remaining sections are organized as follows: Section \ref{section2} presents a detailed description of our model, notations and assumptions. The main results for the CLT of test statistics $U,W,V$ are stated in Section \ref{section4}. 
We also provide numerical studies in Section \ref{section5}.  
Technical proofs are presented in Section \ref{proof of UWV}.  Some derivations and
calculations in Section \ref{proof of UWV}   are postponed to Section \ref{sectionB}. Due to space limitations, simulation results are gathered in the Supplementary Material.

\section{Model}\label{section2}
Throughout the paper, we use bold capital letters and bold italic lowercase letters to represent matrices and vectors, respectively. Scalars are represented by regular letters.  $\boldsymbol{e}_{i}$ denotes a standard basis vector whose components are all zero, except the $i$-th component,
which is equal to 1. We use tr$ (\bbA) $, $ \bbA^{\top} $ and $ \bbA^{\ast} $ to denote the trace, transpose and conjugate transpose of matrix $ \bbA $, respectively. We also use $f'$ to denote the derivative of function $f$, and we use $ \frac{\partial}{\partial z_{1}}f(z_{1},z_{2}) $ to denote the partial derivative of function $ f $ with respect to $ z_{1} $. 
Let $ \left[\bbA \right]_{ij}  $ denote the $ (i,j) $-th entry of the matrix $ \bbA $ and $ \oint_{\mathcal{C}}f(z)dz $ denote the contour integral of $ f(z) $ on the contour $ \mathcal{C} $. Let $ \lambda_{i}^{\bbA} $  be the $ i $th largest eigenvalue of matrix $ \bbA $. 
Weak convergence is denoted by $ \stackrel{d}{\rightarrow}$.  Throughout this paper,  we use $o(1)$ (resp. $o_p (1)$) to denote a negligible scalar (resp. in probability), and the notation $C$ represents a generic constant that may vary from line to line.

\textcolor{black}{In this work, we adopt the notation $ \bbX=(\boldsymbol x_{1},\ldots,\boldsymbol x_{n})=(x_{ij}) $, where $ 1\leq i\leq p $, $ 1\leq j\leq n $.
The singular value decomposition of $\bbT$ is given by
\begin{equation}\label{decT}
	\bbT=\bbV\bbD^{1/2}\bbU^*=
	(\bbV_1,\bbV_2)
	\left( 
	\begin{array}{cc}
		\bbD_{1}^\frac{1}{2} & 0\\
		0 & \bbI_{p-M}
	\end{array}
	\right)
	(\bbU_1,\bbU_2)^{\ast},
\end{equation} 
where
\begin{itemize}
	\item  $\mathbf{V}$ and $\bbU$ are unitary matrices;
	\item $\bbD_{1}$ 
	is a diagonal matrix whose elements  $\al_1\geq\dots\geq\al_K$ are the spiked eigenvalues of $\bSi$ with multiplicities $d_1,\dots,d_K$, respectively. $d_1+d_2+\dots+d_K=M$, and $M$ is a constant. The spiked eigenvalues can be bounded or diverge to infinity.
\end{itemize} }
Then
the corresponding sample covariance matrix $\bbB=\frac{1}{n}\bbT\bbX\bbX^{\ast}\bbT^{\ast}$
is termed the generalized spiked sample covariance matrix. Aligned with the block structure of $\bbD$, we partition $ \bbV=\left(\bbV_{1},\bbV_{2} \right)$,  and $ \bbU=\left(\bbU_{1},\bbU_{2} \right)$, where $\bbV_{1}$  and $\bbU_1$ are $p\times M$ matrices, and define $\bGma=\bbV_{2}\bbD_{2}^{1/2}\bbU_{2}^{\ast} $.
For any matrix $\bbA$ with real eigenvalues, the empirical spectral distribution of $\bbA$ is defined as $$F^{\bbA}(x)=\frac{1}{p}(\text{number of eigenvalues of }  \bbA\leq x).  $$ 
For any function of bounded variation $F$ on the real line, its Stieltjes transform is defined as $$ m_F(z)=\int \frac{1}{\lambda-z}dF(\lambda),~ z\in \mathbb{C}^{+} :=\{z \in \mathbb{C}: \Im z>0\}. $$

The assumptions used to obtain the results in this paper are as follows:
\newtheorem{assumption}{Assumption}[]
{\begin{assumption} \label{ass1}
		$ \{x_{ij},  1\leq i\leq p ,  1\leq j\leq n \} $ {are i.i.d. random variables with common moments} $$ \mathbb{E}x_{ij}=0,\quad  \mathbb{E}\left| x_{ij}\right| ^{2}=1, \quad  \beta_{x}= \mathbb{E}\left|x_{ij}\right| ^{4}- \left|\mathbb{E}x_{ij}^{2}\right|^{2}-2,  \quad \alpha_{x}=\left|\mathbb{E}x_{ij}^{2}\right|^{2}.  $$ 
\end{assumption} }
\begin{assumption} \label{ass2}
	$ \bbT $ is nonrandom, and  $M,~K$ and $d_i (i=1,\dots,K)$ are fixed. As $\min\{p,n\}\to\infty$, the ratio of the dimension-to-sample size $ c_{n}:={p}/{n}\rightarrow c>0. $  $H_n:=F^{\bGma\bGma^{*}}\stackrel{d}{\rightarrow}H$, where $H$ is a distribution function on the real line.
\end{assumption}

Similar to \cite{Silverstein95S} that under Assumptions \ref{ass1} and \ref{ass2}, we have $F^{\bbB}\stackrel{d}{\rightarrow}F^{c,H}$ almost surely, where $F^{c,H}$ is the limiting spectral distribution (LSD) of $\bbB$. 


We first introduce some notations before presenting the main results in the next section. Let $ \underline{F}^{c,H} $ denote the LSD of matrix $ \bbX^{\ast}\bbU_{2}\bbD_{2}\bbU_{2}^{*}\bbX/ n$,  $ c_{nM}=(p-M)/n,~~H_{2n}=F^{\bbD_{2}}. $ 
\begin{table}[]
	\caption{Definitions of the symbology
	}
	\label{notations}
	\centering
	\renewcommand\arraystretch{1.2}
	\begin{tabular}{c|c}
		\hline 
		$ \bbU_{1}=\left(u_{ij} \right)_{i=1,\dots,p;j=1,\dots,M}  $	&  $\mathcal{U}_{i_{1}j_{1}i_{2}j_{2}}=\sum_{t=1}^{p}\overline{u}_{ti_{1}}u_{tj_{1}}u_{ti_{2}}\overline{u}_{tj_{2}}  $\\
		\hline 
		$\phi_n\left(x \right)=x\left(1+c_n\int\frac{t}{x-t}dH_n\left(t \right)  \right)$ & $\phi_{k}=\phi\left(x \right)\mid_{x=\al_{k}}=\al_{k}\left(1+c\int\frac{t}{\al_{k}-t}dH\left(t \right)  \right)$ \\
		\hline 
		$\theta_{k}=\phi_{k}^{2}\underline{m}_{2}\left( \phi_{k}\right)$ & $\nu_{k}=\phi_{k}^{2} \underline{m}^{2}\left(\phi_{k}\right) $  \\
		\hline 
		$ \underline{m}\left( \lambda\right)=\int\frac{1}{x- \lambda }d\underline{F}^{c,H}\left( x\right) $& $ \underline{m}_{2}\left( \lambda\right)=\int\frac{1}{\left( \lambda-x\right) ^{2}}d\underline{F}^{c,H}\left( x\right) $     \\
		\hline
	\end{tabular}
\end{table}	
Moreover, $F^{c_{nM},H_{2n}}$ is the LSD $F^{c,H}$ with $\{c,~H\}$ replaced by $\{c_{nM},~ H_{2n}\}$, and $\mathcal{C}$ is a closed contour in the complex plane enclosing the support of  $ F^{c_n, H_n}$ and it is also enclosed in the analytic area of $ f $. Define	$ s_{k}^{2}= \frac{\left(\al_{x}+1 \right)d_{k} }{\theta_{k}}+\frac{ \beta_{x}\nu_{k} \sum_{j_{1}, j_{2}\in J_{k}}\mathcal{U}_{j_{1}j_{1}j_{2}j_{2}} }{\theta_{k}^{2}}.   $  Other notations used in subsequent sections are defined in Table \ref{notations}. It is worth noting that $ \bbU_{1} $ is the right singular vector matrix of the spiked eigenvalues whose entries $u_{ij}$ are crucial in the CLT established in the paper, and the other symbols in  Table \ref{notations} can be regarded as functions introduced to simplify the following presentation. 



\section{Main results}\label{section4}
In this section, we focus on a hypothesis test that whether $\bSi$ is an identity matrix or follows a generalized spiked model:
\begin{align} \label{hyapp}
	H_{0}:\bSi=\bbI_{p}\quad \text{vs.} \quad
	H_{1}:\bSi=\mathbf{V}\left(\begin{array}{cc}
		\bbD_{1} & 0 \\
		0 & \bbI_{p-M}
	\end{array}\right) \mathbf{V}^{\ast},
\end{align}
where
\begin{itemize}
	\item $M$ is a fixed constant;
	\item $\bbD_1$ is a diagonal matrix of diverging spikes of $\bSi$ ($\al_1\geq\al_2\geq\dots\geq\al_K\rightarrow\infty$);
	\item $\bbV$ is a unitary matrix.
\end{itemize}
In this work, we consider four classical linear hypothesis test statistics  $ U, W, V,R $ and obtain the asymptotic distributions of $ U, W, V,R $ under $ H_{1} $ in (\ref{hyapp}), and we also provide their asymptotic power functions. The specific formulas of $ U, W, V,R $ are given in Section \ref{result UWV}.  
For clarity, we first introduce the following notations: 
\begin{align*}
	&CT(x,c,\tilde{c})= \log( 1+x ) +\frac{-(\sqrt{c}-\frac{1}{\sqrt{c}} )^{2}(\log(1-\sqrt{\tilde{c}c   })+ \sqrt{\tilde{c}c   } )-\sqrt{\tilde{c}}(\sqrt{c}-(\sqrt{c})^{3} )   }{ 1-c},\\
	&I_{1}( f_{U})= \log(1+\varrho(c_{nM})  )-\log(2+c_{nM} )+\sum_{k=1}^{\infty}(\frac{\sqrt{c_{nM}}}{2+c_{nM}} )^{2k}\frac{(2k-1)!}{(k!)^2},\\ 
	&I_{2}( f_{U})=-\sum_{k=1}^{\infty}(\frac{\sqrt{c_{nM}}}{2+c_{nM}} )^{2k}\frac{(2k-1)! }{(k-1)!(k+1)!  },~\varrho(c) =\frac{c+\sqrt{c^{2}+4}}{2},\\
	&J_{1}(f_{U}, f_{U})=(\sum_{k=1}^{\infty}(\frac{\sqrt{c_{nM}}}{2+c_{nM}} )^{2k-1}\frac{(2k-2)! }{k!( k-1) !  } )^{2},~\tilde{c}=\frac{4c}{(2+c+\sqrt{c^{2}+4})^{2} },\\
	&I_{1}( f_{V})=\frac{1}{2+c_{nM}}\sum_{k=0}^{\infty}(\dfrac{\sqrt{c_{nM}}}{2+c_{nM}})^{2k}\frac{(2k)! }{(k!)^2}-\frac{1}{1+\varrho( c_{nM}) }     (\frac{\tilde{c}_{nM}}{-(1-\tilde{c}_{nM})^{2} }+\frac{1}{2(\sqrt{\tilde{c}_{nM}}-1)^{2} }+  \frac{1}{2(\sqrt{\tilde{c}_{nM}}+1)^{2} }   ),\\
	&I_{2}( f_{V})=-\frac{1}{2+c_{nM}}\sum_{k=0}^{\infty}(\frac{\sqrt{c_{nM}}}{2+c_{nM}} )^{2k}\frac{(2k)!}{(k-1)!(k+1)!}, J_{1}(f_{V}, f_{V})=(\frac{1}{2+c_{nM}}\sum_{k=0}^{\infty}(\frac{\sqrt{c_{nM}}}{2+c_{nM}} )^{2k+1} \frac{(2k+1)! }{k!(k+1 )!})^{2}.
\end{align*}
To avoid misunderstandings,  we define the values of $ \varrho\left(c_{n} \right), \varrho\left(c_{nM} \right)   $ to be the same as $ \varrho\left(c \right)  $ above with the substitution of $ c_{n} $ and $ c_{nM} $ for $ c $ in these quantities, respectively. The same substitution also holds for $ \tilde{c}_{n} $.

\subsection{Asymptotic results for test statistics $U,W,V$}

\begin{thm}[U statistics]\label{U}
	Under Assumptions \ref{ass1} and \ref{ass2} with $ c_{n}=p/n\rightarrow c\in\left(0,1 \right)  $,	  we have  under $ H_{1} $ in (\ref{hyapp}),
	\begin{align*}
		\dfrac{U-\left(p-M\right)CT\left(\varrho\left(c_{nM} \right),c_{nM},\tilde{c}_{nM} \right)-\mu_{U} }{\varsigma_{U}}\stackrel{d}{\rightarrow}N\left(0,1 \right), 
	\end{align*}
	where 
	\begin{align*}
		&\mu_{U}=\al_{x}I_{1}\left( f_{U}\right)+\beta_{x} I_{2}\left( f_{U}\right)+\sum_{k=1}^{K}d_{k}\log\left(1+\phi_{n}\left(\al_{k} \right)  \right)+ M\log\left( 1-\sqrt{\tilde{c}_{nM}c_{nM}   }\right), \\
		&\varsigma_{U}^{2}=\sum_{k=1}^{K}\dfrac{ \phi_{n}^{2}\left(\al_{k} \right)}{n\left(1+ \phi_{n}\left(\al_{k} \right)\right)^{2} }s_{k}^{2}+\left( \al_{x}+\beta_{x}+1\right)J_{1}\left(f_{U}, f_{U}\right).
	\end{align*}
\end{thm}

\begin{thm}[W statistics]\label{W}
	Under Assumptions \ref{ass1} and \ref{ass2} with $ c\in\left(0,1 \right)  $,	  we have under $ H_{1} $ in (\ref{hyapp}),
	\begin{align*}
		\dfrac{W-\left(p-M \right)-(\sum_{k=1}^{K}d_{k}\phi_{n}\left(\al_{k} \right) -Mc_{nM}) }{\varsigma_{W}}\stackrel{d}{\rightarrow}N\left(0,1 \right), 
	\end{align*}
	where
	\begin{align*}
		&\varsigma_{W}^{2}=\sum_{k=1}^{K}\dfrac{ \phi_{n}^{2}\left(\al_{k}\right)}{n}s_{k}^{2}+(\al_{x}+\beta_x+1)c_{nM}.	
	\end{align*}
\end{thm}

\begin{thm}[V statistic]\label{V}
	Under Assumptions \ref{ass1} and \ref{ass2} with $  c\in\left(0,1 \right)  $,	  we have  under $ H_{1} $ in (\ref{hyapp}),
	\begin{align*}
		\dfrac{V-\frac{p-M}{1+\varrho\left( c_{nM}\right) }-\mu_{V} }{\varsigma_{V}}\stackrel{d}{\rightarrow}N\left(0,1 \right), 
	\end{align*}
	where 
	\begin{align*}
		\mu_{V}=&\al_{x}I_{1}( f_{V})+\beta_{x} I_{2}( f_{V})+\sum_{k=1}^{K}d_{k}\frac{\phi_{n}(\al_{k} )}{1+\phi_{n}(\al_{k})}-\frac{M( c_{nM}-2) }{2( 1+\varrho(c_{nM}) )(  1-\tilde{c}_{nM} )}-\frac{M}{2}, \\
		\varsigma_{V}^{2}=&\sum_{k=1}^{K}\frac{ \phi_{n}^{2}\left(\al_{k} \right) }{n(1+ \phi_{n}(\al_{k} ))^{4} }s_{k}^{2}+( \al_{x}+\beta_{x}+1)J_{1}(f_{V}, f_{V}).
	\end{align*}
\end{thm}
\begin{remark}
	The proofs of Theorems \ref{U}--\ref{V} are given in Section \ref{proof of UWV}, and to avoid confusion with classical distributions of the Wilks'U test, the Lawley-Hotelling W test, and the Bartlett-Nanda-Pillai V test, we refer the test in Theorems \ref{U}--\ref{V} as `corrected Wilks' likelihood ratio test (CUT)', `corrected Lawley-Hotelling trace test (CWT)' and `corrected Bartlett-Nanda-Pillai trace test (CVT)' instead. From Theorems \ref{U}--\ref{V}, we reject $H_0$ in (\ref{hyapp}) if 
	\begin{align*}
		U>&z_\xi\varsigma_U^0+pCT(\varrho(c_n),c_n,\tilde{c}_n)+\mu_U^0, \\
		W>&z_\xi\sqrt{(\al_{x}+\beta_x+1)c_n}+p,
		\\
		V>&z_\xi\varsigma_{V}^0+\frac{p}{1+\varrho(c_n)}+\mu_V^0,
	\end{align*}  
	where $\xi$ is the significance level of the test and $z_{\xi}$ is the $1-\xi$ quantile of the standard Gaussian distribution $\Phi$. $\varsigma_U^0, \mu_U^0, \varsigma_{V}^0, \mu_V^0$ and the power functions of CUT, CWT, CVT are given in the following Corollaries \ref{Ucase1cor}--\ref{Vcase1cor}.
\end{remark}	
\begin{cor}[Power function of CUT]\label{Ucase1cor}
	Under the same assumptions as in Theorem \ref{U}, we have as $ n\rightarrow \infty $, the power function of CUT $ P_{U}=P(U>z_\xi\varsigma_U^0+pCT(\varrho(c_n),c_n,\tilde{c}_n)+\mu_U^0 )  $ satisfies
	\begin{align*} 
		P_{U}-
		\Phi\left(\dfrac{ \left(p-M \right)CT\left(\varrho\left(c_{nM} \right), c_{nM},\tilde{c}_{nM} \right)-p CT\left(\varrho\left(c_{n} \right), c_{n},\tilde{c}_{n} \right) +A_{1}-z_{\xi}\varsigma_U^0}{\varsigma_{U}}\right)\rightarrow0,
	\end{align*}
	where
	\begin{align*}
		&\mu_{U}^{0}=\al_{x}I_{1}^{0}\left(f_{U} \right)+\beta_{x}I_{2}^{0}\left( f_{U}\right),   \varsigma_{U}^{0}=\sqrt{\left(\al_{x}+\beta_{x}+1 \right)J_{1}^{0}\left(f_{U}, f_{U} \right)},  ~A_{1}=\sum_{k=1}^{K}d_{k}\log\left(1+\phi_{n}\left(\al_{k} \right)  \right)+M\log\left(1-\sqrt{\tilde{c}_{nM}c_{nM}}
		\right),
	\end{align*}
	and $I_{1}^{0}\left( f_{U}\right), I_{2}^{0}\left( f_{U}\right), J_{1}^{0}\left(f_{U}, f_{U}\right) $ are the same as $I_{1}\left( f_{U}\right), I_{2}\left( f_{U}\right), J_{1}\left(f_{U}, f_{U}\right) $ with the substitution of $c_{nM}$ for $c_n$, respectively. 
\end{cor}  

\begin{cor}[Power function of CWT]\label{Wcase1cor}
	Under the same assumptions as in Theorem \ref{W}, we have as $ n\rightarrow\infty, $ the power function of CWT $ P_{W}=P(W>z_\xi\sqrt{(\al_{x}+\beta_x+1)c_n}+p) $ satisfies 
	\begin{align*} 
		P_{W}-\Phi\left(\dfrac{\sum_{k=1}^{K}d_{k}\phi_{n}\left(\al_{k} \right)-Mc_{nM}-M  }{\varsigma_{W} }-z_{\xi}\dfrac{\sqrt{\al_{x}c_{n}+\beta_{x}c_{n}+c_{n}}}{\varsigma_{W}} \right)\rightarrow0,
	\end{align*}
\end{cor}	 
\begin{cor}[Power function of CVT]\label{Vcase1cor}
	Under the same assumptions as in Theorem \ref{V}, we have as $ n\rightarrow\infty, $ the power function of CVT $ P_{V}=P(V>z_\xi\varsigma_{V}^0+\frac{p}{1+\varrho(c_n)}+\mu_V^0) $ satisfies
	\begin{align*} 
		P_{V}-\Phi\left(\frac{ \sum_{k=1}^{K}d_k\frac{\phi_n(\al_k)}{1+\phi_n(\al_k)} +\frac{p-M}{1+\varrho(c_{nM})}- \frac{p}{1+\varrho(c_{n})} -\frac{M}{2}-\frac{\left(c_{nM}-2 \right)M }{2\left(1+\varrho\left(c_{nM} \right)  \right)\left(1-\tilde{c}_{nM} \right)  }  }{\varsigma_{V}}-z_{\xi}\dfrac{\varsigma_{V}^0}{\varsigma_{V}} \right)\rightarrow0,
	\end{align*}
	where  
	\begin{align*}
		&\mu_{V}^{0}=\al_{x}I_{1}^{0}\left(f_{V} \right)+\beta_{x}I_{2}^{0}\left( f_{V}\right),   \varsigma_{V}^{0}=\sqrt{\left(\al_{x}+\beta_{x}+1 \right)J_{1}^{0}\left(f_{V}, f_{V} \right) },  
	\end{align*}
	and $ I_{1}^{0}\left( f_{V}\right), I_{2}^{0}\left( f_{V}\right), J_{1}^{0}\left(f_{V}, f_{V}\right) $ are the same as $I_{1}\left( f_{V}\right), I_{2}\left( f_{V}\right), J_{1}\left(f_{V}, f_{V}\right) $ with the substitution of $c_{nM}$ for $c_n$, respectively. 
\end{cor}  
\begin{remark}
	The proofs of Corollaries \ref{Ucase1cor}--\ref{Vcase1cor} are given in Section \ref{proof of UWV}.
\end{remark}

\subsection{Power analysis}
In this part, we discuss the power functions of $P_U$, $P_W$, $P_V$, and compare them with the power function of Roy's largest root test (RLRT), which we denote it as $P_R$. We assume that $\left\lbrace x_{ij} \right\rbrace $ is real, i.e. $\al_x=1$. The following lemma is borrowed from \cite{DingY18N} and it characterizes the asymptotic distribution of $\lambda_{1}$.
\begin{lemma}[Theorem 2.7 in \cite{DingY18N}]
	Under Assumptions \ref{ass1} and \ref{ass2} and $H_0$ in (\ref{hyapp}), we have
	\begin{align*}
		\frac{\lambda_1-\mu_{r}}{\varsigma_{r}}\stackrel{d}{\rightarrow} F_{TW},
	\end{align*}
	where  $ \mu_r=\left(1+\sqrt{c_n}\right)^{2},  $ $ \varsigma_r=n^{-2/3}\left( 1+\sqrt{c_n}\right)\left(1+\sqrt{c_n^{-1}} \right) ^{1/3} $ and $F_{TW}$ is the Type I Tracy-Widom (TW) distribution.
\end{lemma}
The following lemma characterizes the power function of RLRT, which is given in \cite{Liu23}. Here $t_{\xi}$ is the $1-\xi$ quantile of the TW distribution.
\begin{lemma}[Theorem 4.5 in \cite{Liu23}]
	Under Assumptions
	\ref{ass1} and \ref{ass2}  and $H_1$ in \eqref{hyapp}, if the multiplicity of $\alpha_1$ is one, then the power function of the RLRT $P_R= P(\lambda_1>t_{\xi}{\varsigma_{r}}+\mu_r)$ satisfies 
	\begin{align} 
		P_{R}- \Phi\left( - \frac{t_{\xi}{\varsigma_r}+{\mu_r}-\phi_{n}\left(\al_{1}\right)}{s_1\phi_{n}\left(\al_{1}\right)/\sqrt{n}} \right) \rightarrow 0, \label{powerr}
	\end{align}
	as $n\to\infty$.
\end{lemma}	
To give comparisons of four tests, 
we first define $\varkappa_U,\varkappa_W,\varkappa_V$ and $\varkappa_R$, then  comparisons between $P_U,P_W,P_V,P_R$ equals comparisons between $\varkappa_U,\varkappa_W,\varkappa_V$ and $\varkappa_R$.
We define 
\begin{align*}
	\varkappa_U=&\dfrac{(p-M)CT(\varrho(c_{nM}),c_{nM},\tilde{c}_{nM})-pCT(\varrho(c_{n}),c_{n},\tilde{c}_{n})+A_1   }{   \sqrt{\sum_{k=1}^{K}\frac{\phi_n^2(\al_{k})}{n(1+\phi_n(\al_{k}))^2}s_k^2+(\al_{x}+\beta_{x}+1)J_1(f_U,f_U)} }
	-\dfrac{ z_{\xi}\sqrt{(\al_{x}+\beta_x+1)J_1^0(f_U,f_U)}  }{ \sqrt{\sum_{k=1}^{K}\frac{\phi_n^2(\al_{k})}{n(1+\phi_n(\al_{k}))^2}s_k^2+(\al_{x}+\beta_{x}+1)J_1(f_U,f_U)  }   },\\
	\varkappa_W=&\dfrac{\sum_{k=1}^{K}d_k\phi_{n}(\al_{k})-Mc-M-z_{\xi}\sqrt{(\al_x+\beta_x+1)c_n}   }{ \sqrt{\sum_{k=1}^{K}\frac{\phi_n^2(\al_{k})}{n}s_k^2 +(\al_{x}+\beta_{x}+1)c_{nM}  }  },~ \varkappa_R= \frac{\phi_{n}\left(\al_{1}\right)-{\mu_r}-t_{\xi}{\varsigma_r}}{s_1\phi_{n}\left(\al_{1}\right)/\sqrt{n}},\\
	\varkappa_V=& \dfrac{\sum_{k=1}^{K}d_k\frac{\phi_n(\al_{k})}{1+\phi_n(\al_{k})} +\frac{p-M}{1+\varrho(c_{nM})}-\frac{p}{1+\varrho(c_{n})}-\frac{M}{2} -\frac{(c_n-2)M  }{2(1+\varrho(c_n))(1-\tilde{c}_{n})}  }{ \sqrt{\sum_{k=1}^{K}\frac{\phi_n^2(\al_{k})}{n(1+\phi_n(\al_{k}))^4}s_k^2+(\al_{x}+\beta_{x}+1)J_1(f_V,f_V)  }    } -\frac{z_{\xi}\sqrt{(\al_{x}+\beta_{x}+1)J_1^0(f_V,f_V) }  }{\sqrt{\sum_{k=1}^{K}\frac{\phi_n^2(\al_{k})}{n(1+\phi_n(\al_{k}))^4}s_k^2+(\al_{x}+\beta_{x}+1)J_1(f_V,f_V) }  }.	
\end{align*}
Note that $\{z_{\xi},t_{\xi}{\varsigma_r}$\} are of order $O(1)$, which means of constant order. $\{K,M\}$ are fixed, $0<c<1$,
$ \phi_{n}(\al_{k})=\al_{k}+c+o(1)$ and $s_{k}^{2}= {2d_{k} }+{ \beta_{x}\sum_{j_{1}, j_{2}\in J_{k}}\mathcal{U}_{j_{1}j_{1}j_{2}j_{2}} }+o(1)$. In the sequel, we use the notations  $A_n\simeq B_n$  to denote $A_n= B_n+o(B_n)$.  Then, we have the following conclusions. 
\textcolor{black}{\begin{itemize}
	\item When $M=1$, the divergence rates of 	$\varkappa_U$,   $\varkappa_W$, $\varkappa_V$ and  $\varkappa_R$ are showed in Table \ref{table1}. 
	When $\alpha_1=o(\sqrt{n})$, it is easy to find that  $\varkappa_V<\varkappa_U<\varkappa_W<\varkappa_R$. Moreover, to be noting that when $\alpha_1=\Omega(\sqrt{n})$, $\varkappa_W$ has the same divergence rate as $\varkappa_R$, which means the performance of CWT is as good as RLRT when $\al_1$ is large enough. It is not difficult to prove that $\varkappa_U<\varkappa_W$ even if $\alpha_1=\Omega(\sqrt{n})$.
	Therefore,
	from the formulas of $\varkappa_U$,   $\varkappa_W$, $\varkappa_V$ and  $\varkappa_R$, it is easy to find that when there is one spiked eigenvalue, RLRT has its advantages. As $\al_1$ becomes larger, the advantage of $\varkappa_W$ is highlighted. 
	\begin{table}[t]
		\caption{Divergence rates of  $\varkappa_U$,   $\varkappa_W$, $\varkappa_V$ and  $\varkappa_R$ when $M=1$  }
		\label{table1}
		\centering 
		\begin{tabular}{c|c|c|c|c}
			\hline 
			&{$\varkappa_U$ }&
			{$\varkappa_W$ }& {$\varkappa_V$ }&
			{$\varkappa_R$} \\
			\hline 
			$\alpha_1=o(\sqrt{n})$& $\simeq\log(1+\phi_n(\alpha_1)) $ & $\simeq \phi_n(\alpha_1)$ & $O(1) $ & $\simeq \frac{\sqrt{n}}{s_1} $\\
			\hline 
			$\alpha_1=\Omega(\sqrt{n})$& $\simeq\log(1+\phi_n(\alpha_1)) $ & $\simeq \frac{\sqrt{n}}{s_1}  $ &$ O(1) $ & $ \simeq \frac{\sqrt{n}}{s_1} $\\
			\hline 	
		\end{tabular}
	\end{table}	
		\item  When $M>1$, we take $M=2$ as an illustrative example. The divergence rates of  $\varkappa_U$,   $\varkappa_W$, $\varkappa_V$ and  $\varkappa_R$ when $M=2$ are given in Table \ref{table2}. In addition, we assume two spikes are not equal, and for convenience, we assume they have the same divergence rate, that is, $ \al_2=k_2\al_1$ with some $k_2<1$. We find that, when $\alpha_1=o(\sqrt{n})$, we have $\varkappa_V<\varkappa_U<\varkappa_W<\varkappa_R$. When $\alpha_1=\Omega(\sqrt{n})$, for some suitable value of $k_2$,  $\sqrt{  (s_1^2+k_2^2s_2^2)/ (n(1+k_2^2))} <\sqrt{s_1^2/n}  $ can be satisfied, which means CWT could have higher asymptotic power than RLRT in some cases. CUT and CVT have lower asymptotic powers than RLRT this can also be explained intuitively in the following. 		
		 Actually, to test hypothesis (\ref{hyapp}), CUT and CVT have worse performances than CWT and RLRT, which do not depend on the number of spikes. This is due to the test functions $f_U$ and $f_V$. When $\phi_n(\al_{k})f^{\prime}( \phi_n(\al_{k})) =o(\sqrt{n})$, then asymptotic variance is mainly decided by bulk part. Since $f_U= \log(1+x)$, and $ \phi_n(\al_{k})f_U^{\prime}( \phi_n(\al_{k}))= \phi_n(\al_{k})/(1+ \phi_n(\al_{k}))\leq 1 $, then the value of $\varsigma_{U}^2$ is determined by bulk eigenvalues, which is of constant order, therefore the highest divergence rate of $\varkappa_U$ is $\log(1+\phi_n(\al_k))$; Since  $f_V= x/(1+x)$, and $ \phi_n(\al_{k})f_V^{\prime}( \phi_n(\al_{k}))= \phi_n(\al_{k})/(1+ \phi_n(\al_{k}))^2\leq 1 $, then the value of $\varsigma_{V}^2$ is also determined by bulk eigenvalues,  therefore $\varsigma_V^2$ is of constant order. Therefore the highest divergence rate of $\varkappa_V$ is also of constant order. These analyses also provide explanations for divergence rates in Tables \ref{table1} and \ref{table2}.
	\begin{table}[t]
		\caption{Divergence rates of  $\varkappa_U$,   $\varkappa_W$, $\varkappa_V$ and  $\varkappa_R$ when $M=2$ and $\al_2=k_2\al_1$  }
		\label{table2}
		\centering
		\begin{tabular}{c|c|c|c|c}
			\hline 
			&{$\varkappa_U$ }&
			{$\varkappa_W$ }& {$\varkappa_V$ }&
			{$\varkappa_R$} \\
			\hline 
			$\alpha_1=o(\sqrt{n})$& $\simeq\log(1+\phi_n(\alpha_1)) $ & $\simeq \phi_n(\alpha_1)$ & $O(1) $ & $\simeq \frac{\sqrt{n}}{s_1} $\\
			\hline 
			$\alpha_1=\Omega(\sqrt{n})$& $\simeq\log(1+\phi_n(\alpha_1)) $ & $\simeq \frac{1}{\sqrt{  (s_1^2+k_2^2s_2^2)/ (n(1+k_2^2))}}     $ &$ O(1) $ & $ \simeq \frac{\sqrt{n}}{s_1} $\\
			\hline 
		\end{tabular}
	\end{table}		
\end{itemize}
The following result is a direct consequence of the above analyses, and it holds when the number of spikes $M$ is finite.
\begin{thm}\label{general analysis}
	For four tests CUT, CWT, CVT and RLRT, to detect hypothesis (\ref{hyapp}), when  $\al_1=o(\sqrt{n})$, R has the highest asymptotic power. When  $\al_1=\Omega(\sqrt{n})$, CWT has the highest asymptotic power. CUT and CVT have lower asymptotic powers to detect hypotheses (\ref{hyapp}).
\end{thm}
}
\section{Numerical studies}\label{section5}
In this section, to demonstrate the effectiveness of the proposed CLTs, we provide some short numerical studies.
We examine the following two different distributions of $x_{ij}$:
\begin{itemize}
	\item[$ Dt_1$:] $\left\lbrace x_{ij}\right\rbrace  $ are i.i.d. samples from a standard Gaussian population. 
	\item[$ Dt_2$:]$\left\lbrace x_{ij}\right\rbrace  $ are i.i.d. samples from  $ Gamma(4,0.5)-2 $.
\end{itemize}
In above settings, $ \beta_{x}=0, \frac{3}{2} $ respectively.

In the current numerical studies, the null hypothesis is defined as 
$ H_0: \bSi=\bbI_{p}. $
For the alternative hypothesis, we adopt the following four population covariance matrix structures:

\begin{itemize}
	\item[$H_1$:] Assume that $ \bSi=\bgL_{1}=diag(1+n,\underbrace{1,1,\dots,1}_{p-1}) $ (Model 1).
	\item[$H_2$:] Assume that $\bSi=\bgL_{2}=diag(1+n,1+0.8n,\underbrace{1,1,\dots,1}_{p-2})$ (Model 2).
	\item[$H_3$:] Assume that
	$\bSi=\bbU_{0}\bgL_1\bbU_{0}^{*}, $ 
	where $ \bbU_{0} $ is the left singular vectors of a $ p\times p $ random matrix  with i.i.d. $ N(0,1)$
	entries (Model 3).
	\item[$H_4$:] Assume that  $\bSi=\bbU_{0}\bgL_2\bbU_{0}^{*}, $ 
	and $ \bbU_{0} $ is defined in $ H_{3} $ (Model 4).
\end{itemize}
Due to space limitations, the simulation results are 
gathered in the Supplementary Material.

\section{Technical proofs}\label{proof of UWV}
In this section, we present the proofs of Theorems \ref{U}--\ref{V}, Corollaries \ref{Ucase1cor}--\ref{Vcase1cor}. Before the proofs, some notations and preliminary results are needed. Notations which will be used in the sequel proofs are provided in the following Table \ref{notations2}.
\subsection{Preliminary results}
\begin{table}[t]
	\caption{Definitions of the symbology
	}
	\label{notations2}
	\centering
	\renewcommand\arraystretch{2}
	\begin{tabular}{c|c|c}
		\hline $ m_n(z)=\frac{1}{p}\mathrm{tr}\left(\bbB-z\bbI_{p} \right)^{-1} $&$ \underline{m}_{n}(z)=-\frac{1-c_{n}}{z}+c_{n}m_{n}(z) $
		& $ m_{2n}(z)=\frac{1}{p-M}\mathrm{tr}(\bbS_{22}-z\bbI_{p-M})^{-1}  $\\ \hline$ \underline{m}_{2n}(z)=-\frac{1-c_{nM}}{z}+c_{nM}m_{2n}(z) $ &
		 $ m_{1n0}(z)=\int \frac{1}{x-z}dF^{c_n,H_n}(x) $& $ \underline{m}_{1n0}(z)=-\frac{1-c_{n}}{z}+c_{n}m_{1n0}(z) $\\
		\hline $ m_{2n0}(z)=\int \frac{1}{x-z}dF^{c_{nM},H_{2n}}(x) $&$ \underline{m}_{2n0}(z)=-\frac{1-c_{nM}}{z}+c_{nM}m_{2n0}(z) $&
			$ \bbP_{n}(z)=\left( (1-c_{nM})\bGma\bGma^{*}-zc_{nM}m_{2n0}(z)\bGma\bGma^{*}-z\bbI_{p}\right)^{-1} $ \\
			\hline \multicolumn{3}{c}{
				$ \Theta_{0,n}(z_{1},z_{2})=\frac{\underline{m}_{2n0}^{\prime}(z_{1}) \underline{m}_{2n0}^{\prime}(z_{2})  }{(\underline{m}_{2n0}(z_{1})-\underline{m}_{2n0}(z_{2}) )^{2} }-\frac{1}{(z_{1}-z_{2})^{2}} $,
				$ \Theta_{1,n}(z_{1},z_{2})=\frac{\partial}{\partial z_{2}}\{ \frac{\partial \mathcal{A}_{n}(z_{1},z_{2})}{\partial z_{1}}\frac{1}{1-\al_{x}\mathcal{A}_{n}(z_{1},z_{2})} \} $} \\
		\hline \multicolumn{3}{c}{
			$   \Theta_{2,n}(z_{1},z_{2})=\frac{z_{1}^{2}z_{2}^{2}\underline{m}_{2n0}^{\prime}(z_{1}) \underline{m}_{2n0}^{\prime}(z_{2})}{n}\sum_{i=1}^{p}\left[ \bGma^{*}\bbP_{n}^{2}(z_{1})\bGma\right]  _{ii}\left[ \bGma^{*}\bbP_{n}^{2}(z_{2})\bGma\right]  _{ii}, \vartheta_{n}^{2}=\Theta_{0,n}(z_{1},z_{2})+\al_{x}\Theta_{1,n}(z_{1},z_{2})+\beta_{x}\Theta_{2,n}(z_{1},z_{2})  $ }\\
		\hline \multicolumn{3}{c}{
			$ {\mathcal{A}_{n}(z_{1},z_{2})}=\frac{z_{1}z_{2}}{n}\underline{m}_{2n0}(z_{1}) \underline{m}_{2n0}(z_{2})\mathrm{tr}{\bGma^{*}\bbP_{n}(z_{1})\bGma\bGma^{\top}\bbP_{n}(z_{2})^{\top} \bar{\bGma}}  $}\\
		\hline \multicolumn{3}{c}{$ \makecell[c]{ \mu_{1}=-\frac{\alpha_{x}}{2 \pi i}\cdot\oint_{\mathcal{C}}\frac{  c_{nM} f_{1}(z)\int \underline{m}_{2n0}^{3}(z)t^{2}\left(1+t \underline{m}_{2n0}(z)\right)^{-3} d H_{2n}(t)}{(1-c_{nM} \int \frac{\underline{m}_{2n0}^{2}(z) t^{2}}{\left(1+t \underline{m}_{2n0}(z)\right)^{2}} d H_{2n}(t))(1-\alpha_{x} c_{nM} \int \frac{\underline{m}_{2n0}^{2}(z) t^{2}}{\left(1+t \underline{m}_{2n0}(z)\right)^{2}} d H_{2n}(t)) }dz 
				-\frac{\beta_{x}}{2 \pi i} \cdot \oint_{\mathcal{C}} \frac{c_{nM}  f_{1}(z)\int \underline{m}_{2n0}^{3}(z) t^{2}\left(1+t \underline{m}_{2n0}(z)\right)^{-3} d H_{2n}(t)}{1-c_{nM} \int \underline{m}^{2}_{2n0}(z) t^{2}\left(1+t \underline{m}_{2n0}(z)\right)^{-2} d H_{2n}(t)} dz}$  }\\
		\hline
	\end{tabular}
\end{table}	
Note that $$\sum_{j=1}^{p}f_1\left( \bgl_{j}\right)=p\int f_1\left(x \right)dF^{\bbB}(x),$$ 
and then we define the normalized LSSs as
\begin{align*}
	Y_{1}&= \int f_1\left(x \right)dG_{n}\left( x\right)-\sum_{k=1}^{K}d_kf_1(\phi_n(\al_k))- \frac{1}{2\pi i}\sum_{k=1}^{K}\oint_{\mathcal{C}} f_1(z)\frac{\underline{m}_{2n0}^{\prime}(z) }{1/\alpha_k+\underline{m}_{2n0}(z)}dz,
\end{align*} 
where $$ G_{n}\left( x\right)=p[F^{\bbB}\left(x \right)-F^{c_{n},H_{n}}\left(x \right)] .$$
\begin{assumption}\label{ass4}
	Test function $ f_{1} $ is analytic on a connected open region of the complex plane containing the support of $ F^{c_{n},H_{n}}$ for almost all $n$, where $F^{c_{n},H_{n}} $ to be the same as $F^{c,H} $ above with the subsitution of $c_n$ and $H_n$ for $c$ and $H$. Moreover, we suppose that $$\lim_{\{x_{n},y_n\}\to\infty\atop {x_{n}}/{y_{n}}\rightarrow 1}\frac{f_{1}'\left(x_{n} \right) }{f_{1}'\left(y_{n} \right)}= 1 .$$
\end{assumption}    
In the following theorem, we provide a general CLT result for LSS with test function $f_1$ that satisfies Assumption \ref{ass4}. The population covariance matrix has the structure (\ref{ds}).
\begin{thm}\label{thm1}
Under Assumptions
\ref{ass1}--\ref{ass4}, define $ \varpi_{nk1}^{2}=\frac{\phi_{n}\left(\al_{k} \right)}{\sqrt{n}} f_{1}'\left(\phi_{n}\left(\al_{k} \right)  \right), $ then we have $$ \frac{Y_{1}-\mu_{1}}{\varsigma_{1}}\stackrel{d}{\rightarrow}N\left(0, 1\right),    $$  where
{\begin{align}
	\varsigma_{1}^{2}=\sum_{k=1}^{K}\varpi_{nk1}^{2}s_{k}^{2}  -\frac{1}{4\pi^{2}}\oint_{\mathcal{C}_{1}}\oint_{\mathcal{C}_{2}}f_{1}\left(z_{1} \right)f_{1}\left(z_{2} \right)\vartheta_{n}^{2}dz_{1}dz_{2},\label{thm3.1cov}
\end{align}
and $\mu_1, s_k^2, \vartheta_{n}^2$ are defined in Table \ref{notations2}. $\mathcal{C}_{1}$ and $\mathcal{C}_{2}$ are nonoverlapping and closed contours in the complex plane enclosing the support of  $ F^{c_{n}, H_{n}}$.  $\mathcal{C}_{1}$ and $\mathcal{C}_{2}$ are also enclosed in the analytic area of $ f_1. $  }
\end{thm}	
\begin{remark}\label{remark3.1}
To be noticed that, when all the spiked eigenvalues $\{\alpha_k\}_{1\leq k\leq K}$ tend to infinity, Theorem \ref{thm1} reduces to Theorem 3.1 in \cite{Liu23}.
\end{remark}
\begin{remark}\label{remarkUV}
If we set $f_1$ in Theorem \ref{thm1} to specific functions $f_U$, $f_W$, and $f_V$, and set $\bbD_2$ equals an identity matrix, then Theorems \ref{U}--\ref{V} follow. To guarantee the consistency of the paper, we postpone the proof of Theorem \ref{thm1} to Section \ref{thmf}. 
\end{remark}


\subsection{Proof of Theorem \ref{U}}
Now we prove Theorem \ref{U}. Recall that $$ G_{n}\left( x\right)=p\left[F^{\bbB}\left(x \right)-F^{c_{n},H_{n}}\left(x \right)   \right],  ~ Y_1= \int f_{U}\left(x \right)dG_{n}\left( x\right)-\sum_{k=1}^{K}d_{k}f_{U}\left(\phi_{n}\left(\al_{k} \right)  \right)-\frac{M}{2\pi i}\oint_{\mathcal C}f_{U}\left(z \right)\frac{\underline{m}_{2n0}'(z)}{\underline{m}_{2n0}(z)}dz.$$
When $ f_{U}(x)=\log(1+x), $ after some calculations, we obtain
\begin{gather}
\nonumber\int f_{U}\left(x \right)dG_{n}\left( x\right)=U-p\int  f_{U}(x)dF^{c_{n},H_{n}}(x)\nonumber=U-(p-M)\int f_{U}(x)dF^{c_{nM},H_{2n}}(x),\\
\int f_{U}(x)dF^{c_{nM},H_{2n}}(x)=CT(\varrho(c_{nM}),c_{nM},\tilde{c}_{nM}) \label{Ucenter},\\	
\nonumber\sum_{k=1}^{K}d_{k}f_{U}\left(\phi_{n}\left(\al_{k} \right)  \right)=	\nonumber\sum_{k=1}^{K}d_{k}\log\left(1+\phi_{n}\left(\al_{k} \right)  \right) ,\\
\frac{M}{2\pi i}\oint_{\mathcal C}f_{U}\left(z \right)\frac{\underline{m}_{2n0}'(z)}{\underline{m}_{2n0}(z)}dz=M\log\left(1-\sqrt{\tilde{c}_{nM}c_{nM}} \right). \label{Uextra}
\end{gather} 
For consistency, we postpone the proofs of (\ref{Ucenter}) and (\ref{Uextra}) to Section \ref{sectionB}.
According to Theorem \textcolor{black}{\ref{thm1}} and Theorem A.1 in \cite{wang2013sphericity}, when $ f_{U}(x)=\log\left(1+x \right)  $, we have
\begin{align*}
\frac{U-(p-M)\int f_{U}(x)dF^{c_{nM},H_{2n}}(x)-\mu_{U}}{ \varsigma_{U}}\stackrel{d}{\rightarrow}N(0,1),
\end{align*}	where 
\begin{align}
&\nonumber \mu_{U}=\al_{x}I_{1}\left( f_{U}\right)+\beta_{x} I_{2}\left( f_{U}\right)+\sum_{k=1}^{K}d_{k}\log\left(1+\phi_n(\al_{k}) \right)+M\log( 1-\sqrt{\tilde{c}_{nM}c_{nM}   }), \\
&\nonumber \varsigma_{U}^{2}=\sum_{k=1}^{K}\dfrac{\phi_{n}^{2}(\al_k)}{n(1+\phi_n(\al_k))^2 }s_{k}^{2}+\left( \al_{x}+\beta_{x}+1\right)J_{1}\left(f_{U}, f_{U}\right),  \\
&I_{1}\left( f_{U}\right)=\log\left(1+\varrho(c_{nM}) \right)-\log\left(2+c_{nM} \right)+\sum_{k=1}^{\infty}\left(\frac{\sqrt{c_{nM}}}{2+c_{nM}} \right)^{2k}\dfrac{(2k-1)!}{(k!)^2}\label{Ui1},\\
&I_{2}\left( f_{U}\right)=-\sum_{k=1}^{\infty}\left(\frac{\sqrt{c_{nM}}}{2+c_{nM}} \right)^{2k}\dfrac{\left(2k-1 \right)! }{\left(k-1 \right)!\left(k+1 \right)!  }\label{Ui2},\\
&J_{1}\left(f_{U}, f_{U}\right)=\left(\sum_{k=1}^{\infty}\left(\frac{\sqrt{c_{nM}}}{2+c_{nM}} \right)^{2k-1} \dfrac{\left(2k-2 \right)! }{k!\left(k-1 \right)!}\right)^{2}\label{Uj1}.
\end{align}
Here for consistency, we postpone the proofs of (\ref{Ui1})--(\ref{Uj1}) to Section \ref{sectionB}, and therefore the proof is finished.

\subsection{Proof of Corollary \ref{Ucase1cor}} As the normalized $ U $ statistic tends to a standard normal distribution under $ H_{0} $, that is,
\begin{align*}
\frac{U-p\int f_{U}\left( x\right)dF^{c_{n},H_{n}}-\mu_{U}^{0} }{\varsigma_{U}^{0}}\stackrel{d}{\rightarrow}N(0,1),
\end{align*}	
where  
\begin{align*}
&p\int f_{U}\left( x\right)dF^{c_{n},H_{n}}=pCT\left( \varrho\left( c_{n}\right), c_{n}, \tilde{c}_{n}    \right),~ \mu_{U}^{0}=\al_{x}I_{1}^{0}\left(f_{U} \right)+\beta_{x}I_{2}^{0}\left( f_{U}\right), ~  \varsigma_{U}^{0}=\sqrt{\left(\al_{x}+\beta_{x}+1 \right)J_{1}^{0}\left(f_{U}, f_{U} \right) },  \\
&I_{1}^{0}\left( f_{U}\right)=\log\left(1+\varrho(c_n)\right)-\log\left(2+c_{n} \right)+\sum_{k=1}^{\infty}\left(\frac{\sqrt{c_{n}}}{2+c_{n}} \right)^{2k}\dfrac{(2k-1)!}{k!k!},\\
&I_{2}^{0}\left( f_{U}\right)=-\sum_{k=1}^{\infty}\left(\frac{\sqrt{c_{n}}}{2+c_{n}} \right)^{2k}\dfrac{\left(2k-1 \right)! }{\left(k-1 \right)!\left(k+1 \right)!  },
~J_{1}^{0}\left(f_{U}, f_{U}\right)=\left(\sum_{k=1}^{\infty}\left(\frac{\sqrt{c_{n}}}{2+c_{n}} \right)^{2k-1} \dfrac{\left(2k-2 \right)! }{k!\left(k-1 \right)!}\right)^{2}. 
\end{align*}	
Then we can obtain that $ \xi=P_{H_{0}}\left( U >w \right)=P_{H_{0}}\left(\frac{U -U_{0} }{\varsigma_{U}^{0}}>\frac{w-U_{0}}{\varsigma_{U}^{0}}\right)  $, where $ U_{0}=p\int f_{U}\left( x\right)dF^{c_{n},H_{n}}+ \mu_{U}^{0}$, then critical value $ w=\varsigma_{U}^{0}z_{\xi}+U_{0}. $  Define $ U_{1}=\left( p-M\right)\int f_{U}( x)dF^{c_{nM},H_{2n}}+\mu_{U}  $, then combined with Theorem \ref{U}, we have that the power of test CUT to detect $ H_{1} $ equals 
$
P_U=P_{H_{1}}\left(U>w\right) = P_{H_{1}}\left( \frac{U -U_{1}}{\varsigma_{U}}>\frac{w-U_{1}}{\varsigma_{U}}\right).
$
Since $ \frac{U -U_{1}}{\varsigma_{U}} $ is asymptotically normal distributed, then $ P_U $ is approximate to 
$\Phi\left( \frac{U_{1}-U_{0}-\varsigma_{U}^{0}z_{\xi}  }{\varsigma_{U}}\right).$
Since  $( p-M)\int f_{U}(x)dF^{c_{nM},H_{2n}}- p\int f_{U}\left( x\right)dF^{c_{n},H_{n}}=
\left( p-M\right) CT\left( \varrho\left( c_{nM}\right), c_{nM}, \tilde{c}_{nM}    \right)-pCT\left( \varrho\left( c_{n}\right), c_{n}, \tilde{c}_{n}    \right) $, $\mu_{U}-\mu_{U}^{0}=\al_{x}[  I_{1}\left( f_{U}\right)-  I_{1}^{0}\left( f_{U}\right)] +\beta_{x} [ I_{2}\left( f_{U}\right)-  I_{2}^{0}\left( f_{U}\right)] +
\sum_{k=1}^{K}d_k\log(1+\phi_{n}\left( \al_{k}\right)  )+M\log(1-\sqrt{\tilde{c}_{nM}c_{nM}}) $
and $ I_{1}\left( f_{U}\right)-  I_{1}^{0}\left( f_{U}\right), $ $ I_{2}\left( f_{U}\right)-  I_{2}^{0}\left( f_{U}\right)  $ tend to $ 0 $ as $ n $ tends to infinity, then the proof is finished.
\subsection{Proof of Theorem \ref{W}} Since 
$$  Y_1= \int f_{W}\left(x \right)dG_{n}\left( x\right)-\sum_{k=1}^{K}d_{k}f_{W}\left(\phi_{n}\left(\al_{k} \right)  \right)-\frac{M}{2\pi i}\oint_{\mathcal C}f_{W}\left(z \right)\frac{\underline{m}_{2n0}'(z)}{\underline{m}_{2n0}(z)}dz,  $$ 
when $ f_{W}\left(x \right)=x, $ we obtain
\begin{align}
&\nonumber p\int f_{W}\left(x \right)dF^{c_{n},H_{n}}\left( x\right)=(p-M)\int f_{W}\left(x \right)dF^{c_{nM},H_{2n}}\left( x\right)=p-M, \\
&\sum_{k=1}^{K}d_{k}f_{W}\left(\phi_{n}\left( \al_{k}\right)  \right) = \sum_{k=1}^{K}d_{k} \phi_{n}\left( \al_{k}\right), \label{Wspike} \\
&\frac{M}{2\pi i}\oint_{\mathcal C}f_{W}\left(z \right)\frac{\underline{m}_{2n0}'(z)}{\underline{m}_{2n0}(z)}dz=-Mc_{nM}. \label{Wextra}
\end{align}	
For consistency, we postpone the proof of (\ref{Wextra}) to Section \ref{sectionB}. According to Theorem \ref{thm1}, when $ f_{W}(x)=x $, we have 
\begin{align*}
\dfrac{W-(p-M)\int f_{W}(x)dF^{c_{nM},H_{2n}}(x)-\mu_{W}}{\varsigma_{W}}\stackrel{d}{\rightarrow}N(0,1),
\end{align*}	
where 
\begin{align*}
\mu_{W}=\sum_{k=1}^{K}d_{k} \phi_{n}\left( \al_{k}\right)-Mc_{nM},~ \varsigma_{W}^{2}=\sum_{k=1}^{K}\dfrac{ \phi_{n}^{2}\left(\al_{k}\right)}{n}s_{k}^{2}+\al_{x}c_{nM}+\beta_{x}c_{nM}+c_{nM}.	
\end{align*}	
where $ \mu_{W} $ and $  \varsigma_{W}^{2}$ can be deduced from \cite{wang2013sphericity}, (\ref{Wspike}), and (\ref{Wextra}), therefore the proof of Theorem \ref{W} is finished.
\subsection{Proof of Corollary \ref{Wcase1cor}} 
As the normalized $ W $ statistic tends to a standard normal distribution under $ H_{0} $, that is,
\begin{align*}
\frac{W-p\int f_{W}\left( x\right)dF^{c_{n},H_{n}} }{\varsigma_{W}^{0}}\stackrel{d}{\rightarrow}N(0,1),
\end{align*}	
where  
\begin{align*}
p\int f_{W}\left( x\right)dF^{c_{n},H_{n}}=p,  ~
\varsigma_{W}^{0}=\sqrt{\left(\al_{x}+\beta_{x}+1 \right)c_n }.
\end{align*}	
Then we can obtain that $ P_{H_{0}}\left( W >w \right)=P_{H_{0}}\left(\frac{W -W_{0} }{\varsigma_{W}^{0}}>\frac{w-W_{0}}{\varsigma_{W}^{0}}\right) =\xi  $, where $ W_{0}=p\int f_{W}\left( x\right)dF^{c_{n},H_{n}}$, then critical value $ w=\varsigma_{W}^{0}z_{\xi}+W_{0}. $ Define $ W_{1}=\left( p-M\right)\int f_{W}\left( x\right)dF^{c_{nM},H_{2n}} +\mu_{W}$, then combined with Theorem \ref{W}, we have that the power of test CWT to detect $ H_{1} $ equals 
$
P_W=P_{H_{1}}\left(W >w\right) = P_{H_{1}}\left( \frac{W -W_{1}}{\varsigma_{W}}>\frac{w-W_{1}}{\varsigma_{W}}\right).
$
Since $ \frac{W -W_{1}}{\varsigma_{W}} $ is asymptotically normal distributed, then $ P_W $ is approximate to 	 
$\Phi\left( \frac{W_{1}-W_{0}-\varsigma_{W}^{0}z_{\xi}  }{\varsigma_{W}}\right).$ 
Since  $ ( p-M)\int f_{W}(x)dF^{c_{nM},H_{2n}}- p\int f_{W}( x)dF^{c_{n},H_{n}}=-M $, $\mu_{W}=\sum_{k=1}^{K}d_{k}\phi_{n}\left( \al_{k}\right)-Mc_{nM}, $ therefore $ W_{1}-W_{0}=\sum_{k=1}^{K}d_{k}\phi_{n}\left( \al_{k}\right)-Mc_{nM}-M, $
then the proof is finished.
\subsection{Proof of Theorem \ref{V}}
Since    $$  Y_1= \int f_{V}\left(x \right)dG_{n}\left( x\right)-\sum_{k=1}^{K}d_{k}f_{V}\left(\phi_{n}\left(\al_{k} \right)  \right)-\frac{M}{2\pi i}\oint_{\mathcal C}f_{V}\left(z \right)\frac{\underline{m}_{2n0}'(z)}{\underline{m}_{2n0}(z)}dz, $$ 
when $ f_{V}\left(x \right)=\frac{x}{1+x}, $ we obtain 
\begin{align}
p\int f_{V}\left(x \right)dF^{c_{n},H_{n}}\left( x\right)&=(p-M)\int f_{V}\left(x \right)dF^{c_{nM},H_{2n}}\left( x\right)=\frac{p-M}{1+\varrho\left( c_{nM}\right) }, \label{Vcenter} \\
\nonumber\sum_{k=1}^{K}d_{k}f_{V}\left(\phi_{n}\left( \al_{k}\right)  \right) &= \sum_{k=1}^{K}d_{k} \frac{\phi_{n}\left( \al_{k}\right)}{1+\phi_{n}\left( \al_{k}\right)}, \\
\frac{M}{2\pi i}\oint_{\mathcal C}f_{V}\left(z \right)\frac{\underline{m}_{2n0}'(z)}{\underline{m}_{2n0}(z)}dz&=- \frac{M\left( c_{nM}-2\right) }{2\left( 1+ \varrho\left( c_{nM}\right)\right)\left(1-\tilde{c}_{nM} \right)    }-\frac{M}{2}.\label{Vextra}
\end{align}		
For consistency, we postpone the proofs of (\ref{Vcenter}) and (\ref{Vextra}) to Section \ref{sectionB}. According to Theorem \ref{thm1},   we have
\begin{align*}
\dfrac{V-(p-M)\int f_{V}(x)dF^{c_{nM},H_{2n}}(x)-\mu_{V}}{\varsigma_{V}}\stackrel{d}{\rightarrow}N(0,1),
\end{align*}	
where 
\begin{align}
&\nonumber\mu_{V}=\al_{x}I_{1}(f_{V})+\beta_{x} I_{2}( f_{V})+\sum_{k=1}^{K}d_{k}f_V(\phi_{n}(\al_{k} ))-\frac{M( c_{nM}-2) }{2(1+\varrho(c_{nM}))(1-\tilde{c}_{nM})}-\frac{M}{2}, \\
&\nonumber\varsigma_{V}^{2}=\sum_{k=1}^{K}\dfrac{ \phi_{n}^{2}\left(\al_{k} \right) }{n\left(1+ \phi_{n}\left(\al_{k} \right)\right)^{4} }s_{k}^{2 }+\left( \al_{x}+\beta_{x}+1\right)J_{1}\left(f_{V}, f_{V}\right),  \\
&I_{1}( f_{V})=\frac{1}{2+c_{nM}}\sum_{k=0}^{\infty}(\frac{\sqrt{c_{nM}}}{2+c_{nM}} )^{2k}\frac{\left(2k \right)! }{(k!)^2}-
\frac{1}{1+\varrho(c_{nM}) }     (\frac{\tilde{c}_{nM}}{-(1-\tilde{c}_{nM} )^{2} }+\frac{1}{2(\sqrt{\tilde{c}_{nM}}-1)^{2} }+  \frac{1}{2(\sqrt{\tilde{c}_{nM}}+1 )^{2} }),\label{Vi1}\\
&I_{2}( f_{V})=-\frac{1}{2+c_{nM}}\sum_{k=0}^{\infty}(\frac{\sqrt{c_{nM}}}{2+c_{nM}} )^{2k}\frac{(2k)!}{(k-1)!(k+1)!}\label{Vi2},\\
&J_{1}\left(f_{V}, f_{V}\right)=(\frac{1}{2+c_{nM}}\sum_{k=0}^{\infty}(\frac{\sqrt{c_{nM}}}{2+c_{nM}} )^{2k+1} \frac{(2k+1)! }{k!(k+1)!})^{2}.\label{Vj1}
\end{align}	
We postpone the proofs of (\ref{Vi1})- -(\ref{Vj1}) to Section \ref{sectionB}, then the proof is finished.

\subsection{Proof of Corollary \ref{Vcase1cor}}
As the normalized $ V $ statistic tends to a standard normal distribution under $ H_{0} $, that is,
\begin{align*}
\frac{V-p\int f_{V}\left( x\right)dF^{c_{n},H_{n}}-\mu_{V}^{0} }{\varsigma_{V}^{0}}\stackrel{d}{\rightarrow}N(0,1),
\end{align*}	
where  
\begin{align*}
&p\int f_{V}\left( x\right)dF^{c_{n},H_{n}}=\frac{p}{1+\varrho\left(c_{n} \right)}, ~\mu_{V}^{0}=\al_{x}I_{1}^{0}\left(f_{V} \right)+\beta_{x}I_{2}^{0}\left( f_{V}\right), ~ \varsigma_{V}^{0}=\sqrt{\left(\al_{x}+\beta_{x}+1 \right)J_{1}^{0}\left(f_{V}, f_{V} \right)},  \\
&I_{1}^{0}( f_{V})=\frac{1}{2+c_{n}}\sum_{k=0}^{\infty}(\frac{\sqrt{c_{n}}}{2+c_{n}} )^{2k}\frac{\left(2k \right)! }{(k!)^2}-\frac{1}{1+\varrho( c_{n}) }(\frac{\tilde{c}_{n}}{-(1-\tilde{c}_{n})^{2} }+\frac{1}{2(\sqrt{\tilde{c}_{n}}-1)^{2} }+  \frac{1}{2(\sqrt{\tilde{c}_{n}}+1 )^{2} })  ,\\
&I_{2}^{0}( f_{V})=-\frac{1}{2+c_{n}}\sum_{k=0}^{\infty}(\frac{\sqrt{c_{n}}}{2+c_{n}} )^{2k}\frac{(2k)!}{(k-1)!(k+1)!}, ~J_{1}^{0}\left(f_{V}, f_{V}\right)=\left(\dfrac{1}{2+c_{n}}\sum_{k=0}^{\infty}\left(\frac{\sqrt{c_{n}}}{2+c_{n}} \right)^{2k+1} \dfrac{\left(2k+1 \right)! }{k!\left(k+1 \right)!}\right)^{2}. 
\end{align*}	
Then we can obtain that $ P_{H_{0}}\left( V >w \right)=P_{H_{0}}\left(\frac{V -V_{0} }{\varsigma_{V}^{0}}>\frac{w-V_{0}}{\varsigma_{V}^{0}}\right) =\xi  $, where $ V_{0}=p\int f_{V}\left( x\right)dF^{c_{n},H_{n}}+ \mu_{V}^{0}$, then critical value $ w=\varsigma_{V}^{0}z_{\xi}+V_{0}. $ Define $ V_{1}=( p-M)\int f_{V}( x)dF^{c_{nM},H_{2n}}+\mu_{V}  $, then combined with Theorem \ref{V}, we have that the power of test CVT to detect $ H_{1} $ equals 
$
P_V=P_{H_{1}}\left(V >w\right) = P_{H_{1}}\left( \frac{V-V_{1}}{\varsigma_{V}}>\frac{w-V_{1}}{\varsigma_{V}}\right).
$
Since	 $ \frac{V-V_{1}}{\varsigma_{V}} $ is asymptotically normal distributed, then $ P_V $ is approximate to 
$\Phi\left( \frac{V_{1}-V_{0}-\varsigma_{V}^{0}z_{\xi}  }{\varsigma_{V}}\right). $
Since $(p-M)\int f_{V}(x)dF^{c_{nM},H_{2n}}-p\int f_{V}( x)dF^{c_{n},H_{n}}=\frac{p-M}{1+\varrho( c_{nM}) }-\frac{p}{1+\varrho( c_{n}) }
$, $\mu_{V}-\mu_{V}^{0}=\al_{x}\left[  I_{1}\left( f_{V}\right)-  I_{1}^{0}\left( f_{V}\right)\right] +\beta_{x} \left[ I_{2}\left( f_{V}\right)-  I_{2}^{0}\left( f_{V}\right)\right] +\sum_{k=1}^{K}d_{k}\frac{\phi_{n}\left( \al_{k}\right) }{1+\phi_{n}\left( \al_{k}\right)}  - \frac{M\left( c_{nM}-2\right) }{2\left( 1+ \varrho\left( c_{nM}\right)\right)\left(1-\tilde{c}_{nM} \right)    }-\frac{M}{2}$,
and as $ n $ tends to infinity, $ I_{1}\left( f_{V}\right)-  I_{1}^{0}\left( f_{V}\right) $,  $ I_{2}\left( f_{V}\right)-  I_{2}^{0}\left( f_{V}\right)  $ tend to $ 0 $, then $V_1-V_0$ tends to $\sum_{k=1}^{K}d_{k}\frac{\phi_{n}\left( \al_{k}\right) }{1+\phi_{n}\left( \al_{k}\right)} +\frac{p-M}{1+\varrho\left( c_{nM}\right) }-\frac{p}{1+\varrho\left( c_{n}\right) } - \frac{M\left( c_{nM}-2\right) }{2\left( 1+ \varrho\left( c_{nM}\right)\right)\left(1-\tilde{c}_{nM} \right)    }-\frac{M}{2}$,  
then the proof is finished.	
\subsection{Proof of Theorem \ref{thm1}}\label{thmf}
	First, the sample covariance matrix requires block-wise partitioning. For the population covariance matrix $\bSi=\bbT\bbT^{\ast}$, we consider the corresponding sample covariance matrix $ \bbB=\bbT\bbS_{x}\bbT^{\ast} $, where $ \bbS_{x}=\frac{1}{n}\bbX\bbX^{\ast}  $. 
	By singular value decomposition of $ \bbT $ (see \eqref{decT}),
	\begin{equation*}
		\bbB=\bbV\left(\begin{array}{cc}
			\bbD_{1}^\frac{1}{2}\bbU_{1}^{\ast}\bbS_{x}\bbU_{1}\bbD_{1}^\frac{1}{2},  &\bbD_{1}^\frac{1}{2}\bbU_{1}^{\ast}\bbS_{x}\bbU_{2}\bbD_{2}^\frac{1}{2}\\
			\bbD_{2}^\frac{1}{2}\bbU_{2}^{\ast}\bbS_{x}\bbU_{1}\bbD_{1}^\frac{1}{2},	&\bbD_{2}^\frac{1}{2}\bbU_{2}^{\ast}\bbS_{x}\bbU_{2}\bbD_{2}^\frac{1}{2}
		\end{array} \right)\bbV^{\ast}.
	\end{equation*}
	Note that
	\begin{equation*}
		\bbS=\left(\begin{array}{cc}
			\bbD_{1}^\frac{1}{2}\bbU_{1}^{\ast}\bbS_{x}\bbU_{1}\bbD_{1}^\frac{1}{2},  &\bbD_{1}^\frac{1}{2}\bbU_{1}^{\ast}\bbS_{x}\bbU_{2}\bbD_{2}^\frac{1}{2}\\
			\bbD_{2}^\frac{1}{2}\bbU_{2}^{\ast}\bbS_{x}\bbU_{1}\bbD_{1}^\frac{1}{2},	&\bbD_{2}^\frac{1}{2}\bbU_{2}^{\ast}\bbS_{x}\bbU_{2}\bbD_{2}^\frac{1}{2}
		\end{array} \right)	\triangleq\left(\begin{array}{cc}\bbS_{11},&\bbS_{12}\\\bbS_{21},&\bbS_{22}\end{array} \right). 
	\end{equation*}
	where $ \bbB $ and $ \bbS $ have the same eigenvalues. Let $ \tilde{\bgl}_{j} $ be the eigenvalues of $ \bbS_{22} $ so that the LSS of $ \bbS_{22} $ is $ \sum_{j=1}^{p-M}f( \tilde{\bgl}_{j}) $.
	Before we prove Theorem \ref{thm1}, a technical   lemma is needed, which measures an asymptotic difference between $\sum_{j=M+1}^{p}f(\lambda_{j})$ and $\sum_{j=1}^{p-M}f(\tilde{\lambda}_j)$.
	\begin{lemma}\label{lem1} Under Assumptions \ref{ass1} and \ref{ass2}, we have
		\begin{align*}
			\sum_{j=M+1}^{p}f(\lambda_{j})-\sum_{j=1}^{p-M}f(\tilde{\lambda}_j)=\frac{1}{2\pi i}\sum_{i=1}^{M}\oint_{\mathcal{C}} f(z)\frac{\underline{m}_{2n0}^{\prime}(z) }{1/\alpha_i+\underline{m}_{2n0}(z)}dz+o_p(1).
		\end{align*}	
	\end{lemma}
	\begin{remark}
		When spiked eigenvalues $\alpha_i$ tend to infinity, the result above reduces to Lemma 6.2 in \cite{Liu23}. To guarantee the coherence of the paper, we postpone the proof of Lemma \ref{lem1} until after the proof of Theorem \ref{thm1}. 
	\end{remark}
\textcolor{black}{Now, we continue to the proof of Theorem \ref{thm1}. The proof of Theorem \ref{thm1} builds on the decomposition of LSSs. It is worth noting that, one can follow the same lines of Theorem 3.1 in \cite{Liu23} except that the term $ \sum_{j=M+1}^{p}f(\lambda_{j})-\sum_{j=1}^{p-M}f(\tilde{\lambda}_j) $ is replaced by $\frac{1}{2\pi i}\sum_{i=1}^{M}\oint_{\mathcal{C}} f(z)\frac{\underline{m}_{2n0}^{\prime}(z) }{1/\alpha_i+\underline{m}_{2n0}(z)}dz+o_p(1)$. Therefore, we omit the rest of the proof.}

\section{Some deviations and calculations}\label{sectionB}
Some derivations and calculations in Section \ref{proof of UWV}   are postponed to this section. 

\begin{proof}[Proof of Lemma \ref{lem1}:]	
	We denote $L_1=\sum_{j=M+1}^{p}f(\lambda_j)$, and $L_2=\sum_{j=1}^{p-M}f(\tilde{\lambda}_j)$.	By applying the block matrix inversion formula to $m_n$, we can obtain
	\begin{align}\label{L12}
		L_{1}-L_{2}=-\frac{1}{2\pi i}\oint_{\mathcal C}f\left(z \right)\left(T_{1}-T_{2} \right) dz,
	\end{align}
	where 
	\begin{align*}
		T_{1} &=\mathrm{tr}\left( \bbS_{11}-z\bbI_{M}-\bbS_{12}\left( \bbS_{22}-z\bbI_{p-M}\right)^{-1}\bbS_{21} \right)^{-1},\\ 	
		T_{2} &=-\mathrm{tr}\left[\left( \bbS_{11}-z\bbI_{M}-\bbS_{12}\left(\bbS_{22}-z\bbI_{p-M} \right)^{-1}\bbS_{21} \right)^{-1}\bbS_{12}\left(\bbS_{22}-z\bbI_{p-M} \right)^{-2}\bbS_{21}\right],
	\end{align*}
	which, together with the notation $$\bUps_n:=\frac{1}{n}\bbD_{1}^{\frac{1}{2}}\bbU_{1}^{\ast}\bbX(\frac{1}{n}\bbX^{\ast}\bbU_{2}\bbD_{2}\bbU_{2}^{*}\bbX-z\bbI_{n}) ^{-1}  \bbX^{\ast}\bbU_{1}\bbD_{1}^{\frac{1}{2}}, $$ implies that
	\begin{align*}
		T_{1}&=-z^{-1}\mathrm{tr}\left(\bbI_{M}+\bUps_n \right)^{-1},
		~	T_{2}=z^{-1}\mathrm{tr}\left[\left(\bbI_{M}+\bUps_n \right)^{-1}\bbS_{12}\left(\bbS_{22}-z\bbI_{p-M} \right)^{-2}\bbS_{21}\right].
	\end{align*}
	$ \underline{m}_{2n}= \underline{m}_{2n}(z)$ denotes the Stieltjes transform of $ F^{\frac{1}{n}\bbX^{\ast}\bbU_{2}\bbD_{2}\bbU_{2}^{*}\bbX} $. Thus,  we have that $\underline{m}_{2n}(z)-\underline{m}(z)=o_p(1)$ for any $z\in\mathcal{C}$. 
	From Theorem 3.1 of \cite{JiangB21G}, we know that
	\begin{align} 
		\frac{1}{n}\bbU_{1}^{\ast}\bbX(\frac{1}{n}\bbX^{\ast}\bbU_{2}\bbD_{2}\bbU_{2}^{*}\bbX-z\bbI_{n}) ^{-1}\bbX^{\ast}\bbU_{1}
		=\underline{m}_{2n}\left(z \right)\bbI_{M}+O_{p}(n^{-\frac{1}{2}}).
		\label{6}
	\end{align} 
	Then we have
	\begin{align*}
		T_1=&-\frac{1}{z}\mathrm{tr}(\bbD_1^{-1} + \frac{1}{n}\bbU_{1}^{\ast}\bbX(\frac{1}{n}\bbX^{\ast}\bbU_{2}\bbD_{2}\bbU_{2}^{*}\bbX-z\bbI_{n}) ^{-1}  \bbX^{\ast}\bbU_{1}  )^{-1}\bbD_1^{-1}
		= -\frac{1}{z}\sum_{i=1}^{M}\frac{1}{1+\al_i\underline{m}_{2n}(z)}+O_p(\frac{1}{\sqrt{n}}).
	\end{align*}
	For $T_2$, since
	\begin{align*}
		&(\bbD_1^{-1}+\frac{1}{n}\bbU_{1}^{\ast}\bbX(\frac{1}{n}\bbX^{\ast}\bbU_{2}\bbD_{2}\bbU_{2}^{*}\bbX-z\bbI_{n})^{-1}  \bbX^{\ast}\bbU_{1})^{-1}
		=
		\begin{pmatrix}
			\frac{1}{\al_1}+\underline{m}_{2n} &  & \\  & \ddots &  \\
			&  & \frac{1}{\al_M}+\underline{m}_{2n}
		\end{pmatrix}
		^{-1}+O_p(\frac{1}{\sqrt{n}}),
	\end{align*}
	and from \cite{Liu23},
	\begin{align*}
		&\bbD_{1}^{-1/2}\bbS_{12}\left(\bbS_{22}-z\bbI_{p} \right)^{-2}\bbS_{21}\bbD_{1}^{-1/2}
		= cm_{2n}\left(z \right)\bbI_{M}+zcm_{2n}'\left(z \right)\bbI_{M}+O_p(n^{-3/2})
		=\underline{m}_{2n}\left(z \right)\bbI_{M}+z\underline{m}_{2n}'\left(z \right)\bbI_{M}+O_p(n^{-3/2}),  
	\end{align*}
	then we have 
	\begin{align*}
		T_2=&\frac{1}{z}\mathrm{tr}
		\begin{pmatrix}
			\frac{1}{\al_1}+\underline{m}_{2n} &  & \\  & \ddots &  \\
			&  & \frac{1}{\al_M}+\underline{m}_{2n}
		\end{pmatrix}
		^{-1}
		\begin{pmatrix}
			\underline{m}_{2n}+z\underline{m}_{2n}^{\prime} &  & \\  & \ddots &  \\
			&  & \underline{m}_{2n}+z\underline{m}_{2n}^{\prime}
		\end{pmatrix}+O_p(n^{-2})= \frac{1}{z}\sum_{i=1}^{M}\frac{\underline{m}_{2n}+z\underline{m}_{2n}^{\prime}}{\frac{1}{\al_i}+\underline{m}_{2n}}+O_p(n^{-2}).
	\end{align*}
	Therefore
	\begin{align*}
		L_1-L_2=&-\frac{1}{2\pi i}\oint_{\mathcal{C}} f(z)(T_1-T_2)dz\\=&
		\frac{1}{2\pi i}\sum_{i=1}^{M}\oint_{\mathcal{C}} \frac{f(z)}{z}\frac{ 1 }{1+\al_i\underline{m}_{2n}(z)}dz+ \frac{1}{2\pi i}\sum_{i=1}^{M}\oint_{\mathcal{C}} \frac{f(z)}{z}\frac{ \underline{m}_{2n}(z)+z\underline{m}_{2n}^{\prime}(z) }{\frac{1}{\al_i}+\underline{m}_{2n}(z)}dz+O_p(\frac{1}{\sqrt{n}}).
	\end{align*}
	Since 
	\begin{align*}
		&\frac{1}{2\pi i}\sum_{i=1}^{M}\oint_{\mathcal{C}_{2n}} \frac{f(z)}{z}\frac{ 1 }{1+\al_i\underline{m}_{2n}(z)}dz\\
		=&\frac{1}{2\pi i}\sum_{i=1}^{M}\oint_{\mathcal{C}_{2n}} \frac{f(z)}{z}\frac{ 1 }{1+\al_i\underline{m}_{2n0}(z)} \frac{1+\al_i\underline{m}_{2n0}(z)}{1+\al_i\underline{m}_{2n}(z)}dz=
		\frac{1}{2\pi i}\sum_{i=1}^{M}\oint_{\mathcal{C}_{2n0}} \frac{f(z)}{z}\frac{ 1 }{1+\al_i\underline{m}_{2n0}(z)}dz+o_p(1), 
	\end{align*}
	and 
	\begin{align*}
		&\frac{1}{2\pi i}\sum_{i=1}^{M}\oint_{\mathcal{C}_{2n}} \frac{f(z)}{z}\frac{ \underline{m}_{2n}(z)+z\underline{m}_{2n}^{\prime}(z) }{\frac{1}{\al_i}+\underline{m}_{2n}(z)}dz
		= \frac{1}{2\pi i}\sum_{i=1}^{M}\oint_{\mathcal{C}_{2n}} \frac{f(z)}{z}\frac{ \underline{m}_{2n0}(z)+z\underline{m}_{2n0}^{\prime}(z) }{\frac{1}{\al_i}+\underline{m}_{2n0}(z)} \frac{ \frac{1}{\al_i}+\underline{m}_{2n0}(z) }{\frac{1}{\al_i}+\underline{m}_{2n}(z)  } \frac{\underline{m}_{2n}(z)+z\underline{m}_{2n}^{\prime}(z)}{\underline{m}_{2n0}(z)+z\underline{m}_{2n0}^{\prime}(z)}dz\\
		=&\frac{1}{2\pi i}\sum_{i=1}^{M}\oint_{\mathcal{C}_{2n}} \frac{f(z)}{z}\frac{ \underline{m}_{2n0}(z)+z\underline{m}_{2n0}^{\prime}(z) }{\frac{1}{\al_i}+\underline{m}_{2n0}(z)}dz+o_p(1),
	\end{align*}
	then the proof is finished. 
\end{proof}

\begin{proof}[Proof of (\ref{Ucenter}):]	Since $ f_{U}\left(x \right)=\log\left( 1+x\right),   $ then 
	\begin{align*}
		\int f_{U}(x)dF^{c_{n},H_{n}}(x)=\int_{a(c_{n})}^{b(c_{n})}\log\left(1+x \right)\frac{1}{2\pi xc_{n}}\sqrt{\left(b(c_{n})-x \right) \left(x-a(c_{n}) \right)  }dx, 
	\end{align*}
	where $ a\left(c_{n} \right)=\left( 1-\sqrt{c_{n}}\right)^{2}   $, $ b\left(c_{n} \right)=\left( 1+\sqrt{c_{n}}\right)^{2}   $. By using the variable change $ x=1+c_{n}- 2\sqrt{c_{n}} \cos\left(\theta \right)  $, $ 0\leq\theta\leq\pi, $ we have	
	\begin{align*}
		\int f_{U}(x)dF^{c_{n},H_{n}}(x)&=\frac{1}{2\pi c_{n}}\int_{0}^{\pi}\dfrac{\log\left( 2+c_{n}-2\sqrt{c_{n}}\cos\left(\theta\right) \right) }{1+c_{n}-2\sqrt{c_{n}}\cos\left(\theta\right)}4c_{n}\sin^{2}\left(\theta \right)d \theta\\
		&=\frac{1}{2\pi }\int_{0}^{2\pi}\dfrac{2\sin^{2}\left(\theta \right) }{1+c_{n}-2\sqrt{c_{n}}\cos\left(\theta\right)}\log\left( 2+c_{n}-2\sqrt{c_{n}}\cos\left(\theta\right) \right)d \theta.
	\end{align*}
	Do the transformation $ 2+c_{n}-2\sqrt{c_{n}}\cos\left(\theta\right)=\left(1+\varrho\left(c_{n} \right)  \right)(\frac{2+c_{n}}{1+\varrho(c_{n} ) }-\frac{2\sqrt{c_{n}}}{1+\varrho(c_{n} ) }\cos(\theta ) ),   $   and let $ \frac{2+c_{n}}{1+\varrho\left(c_{n} \right) }=1+\tilde{c}_{n}, $ $ \frac{\sqrt{c_{n}}}{1+\varrho\left(c_{n} \right) }=\sqrt{\tilde{c}_{n}}, $ then we obtain $ \varrho\left(c_{n} \right) =\frac{c_{n}+\sqrt{c_{n}^{2}+4}}{2}, $  $ \tilde{c}_{n}=\frac{4c_{n}}{( 2+c_{n}+\sqrt{c_{n}^{2}+4})^{2} }. $
	Therefore $\int f_{U}(x)dF^{c_{n},H_{n}}(x)  $ equals
	\begin{align*}
		&\frac{1}{2\pi }\int_{0}^{2\pi}\dfrac{2\sin^{2}\left(\theta \right)\log\left(1+ \varrho\left(c_{n} \right) \right) }{1+c_{n}-2\sqrt{c_{n}}\cos\left(\theta\right)}d\theta+ \frac{1}{2\pi }\int_{0}^{2\pi}\dfrac{2\sin^{2}\left(\theta \right) }{1+c_{n}-2\sqrt{c_{n}}\cos\left(\theta\right)}\log( \dfrac{2+c_{n}}{1+\varrho(c_{n} ) }-\dfrac{2\sqrt{c_{n}}}{1+\varrho(c_{n} ) }\cos(\theta ) ) d\theta.
	\end{align*}
	For the first integral,
	\begin{align*}
		&\frac{1}{2\pi }\int_{0}^{2\pi}\dfrac{2\sin^{2}\left(\theta \right)\log\left(1+ \varrho\left(c_{n} \right) \right) }{1+c_{n}-2\sqrt{c_{n}}\cos\left(\theta\right)}d\theta
		=\frac{\log\left(1+ \varrho\left(c_{n} \right) \right) }{2\pi }\int_{0}^{2\pi}\dfrac{2\sin^{2}\left(\theta \right) }{1+c_{n}-2\sqrt{c_{n}}\cos\left(\theta\right)}d\theta\\
		=&-\frac{1}{2}\frac{\log\left(1+ \varrho\left(c_{n} \right) \right) }{2\pi i }\oint_{\left| z\right|=1 }\dfrac{\left( z-\frac{1}{z}\right)^{2} }{\left( z-\sqrt{c_{n}}\right)\left(1-\sqrt{c_{n}}z \right)  }dz
		=-\frac{1}{2}\frac{\log\left(1+ \varrho\left(c_{n} \right) \right) }{2\pi i }\oint_{\left| z\right|=1 }\dfrac{ z^{4}-2z^{2}+1}{z^{2}\left( z-\sqrt{c_{n}}\right)\left(1-\sqrt{c_{n}}z \right)  }dz
	\end{align*}	
	Since $ c_{n}<1, $ thus $ \sqrt{c_{n}} $ and $ 0 $ are poles. The residues are $ \frac{1-c_{n}}{c_{n}} $ and $ -\frac{1+c_{n}}{c_{n}} $, respectively. By the residue theorem, we obtain the first integral is $ \log(1+\frac{c_{n}+\sqrt{c_{n}^{2}+4}}{2})  $.
	For the second integral, 
	\begin{align*}
		&\frac{1}{2\pi }\int_{0}^{2\pi}\frac{2\sin^{2}(\theta ) }{1+c_{n}-2\sqrt{c_{n}}\cos(\theta)}\log( \frac{2+c_{n}}{1+\varrho(c_{n}) }-\dfrac{2\sqrt{c_{n}}}{1+\varrho(c_{n}) }\cos(\theta )) d\theta\\
		=&\frac{1}{2\pi }\int_{0}^{2\pi}\frac{2\sin^{2}(\theta) }{1+c_{n}-2\sqrt{c_{n}}\cos(\theta)}\log(1+\tilde{c}_{n}-2\sqrt{\tilde{c}_{n}}\cos(\theta)) d\theta\\
		=&\frac{-(\sqrt{c_{n}}-\frac{1}{\sqrt{c_{n}}} )^{2}(\log(1-\sqrt{ \tilde{c}_{n}c_{n}})+\sqrt{ \tilde{c}_{n}c_{n}} )-\sqrt{\tilde{c}_{n} } ( \sqrt{c_{n}}-\left(\sqrt{c_{n}} \right)^{3} ) }{1-c_{n}},
	\end{align*}
	where the last integral is calculated in \cite{10.1214/aop/1078415845}. Note that all the $ c_{n} $ in the formulas above should be replaced by $ c_{nM} $ since the calculation is on the bulk part of LSS. Collecting the two integrals leads to the desired formula for $ p\int f_{U}(x)dF^{c_{n},H_{n}}(x) $.
\end{proof}
\begin{proof}[Proof of (\ref{Uextra}):]
	First, we consider $ \oint_{\mathcal C}f_{U}\left(z \right)\frac{\underline{m}'(z)}{\underline{m}(z)}dz $, we have
	\begin{align}
		&\nonumber\oint_{\mathcal C}f_{U}\left(z \right)\frac{\underline{m}'(z)}{\underline{m}(z)}dz
		=\nonumber\oint_{\mathcal C}f_{U}\left(z \right)d\log\underline{m}\left(z \right)=\nonumber-\oint_{\mathcal C}f_{U}^{'}\left(z \right)\log\underline{m}\left(z \right)dz \\
		=&\int_{a(c)}^{b(c)}f_{U}^{'}\left(x \right)\left[\log\underline{m}(x+i\varepsilon)-\log\underline{m}(x-i\varepsilon) \right]dx
		=2i 	\int_{a(c)}^{b(c)}f_{U}^{'}\left(x \right) \Im \log\underline{m}(x+i\varepsilon)dx. \label{26}
	\end{align}	
	Here, $a(c)=(1-\sqrt{c})^{2}  $ and $b(c)=(1+\sqrt{c})^{2}  $. Since $\underline{m}\left(z \right)=-\frac{1-c}{z}+cm\left(z \right), $
	and under $ H_{1} $, we have $ \underline{m}\left(z \right)=\frac{-(z+1-c)+\sqrt{(z-1-c)^{2}-4c}}{2z} . $
	As $z \rightarrow x \in$ $[a(c), b(c)]$, we obtain
	$ \underline{m}\left(x \right)= \frac{-(x+1-c)+\sqrt{4c-(x-1-c)^{2}}i}{2x}.$
	Therefore,
	\begin{align*}
		&\int_{a(c)}^{b(c)}f_{U}^{'}\left(x \right) \Im \log\underline{m}(x+i\varepsilon)dx
		=\int_{a(c)}^{b(c)} f_{U}^{'}(x)\tan ^{-1}(\frac{\sqrt{4 c-(x-1-c)^{2}}}{-(x+1-c)}) d x\\
		=&\left.\tan ^{-1}(  \frac{\sqrt{4 c-(x-1-c)^{2}}}{-(x+1-c)})  f_{U}(x)\right|_{a(c)} ^{b(c)}-\int_{a(c)}^{b(c)} f_{U}(x) d \tan ^{-1}(\frac{\sqrt{4 c-(x-1-c)^{2}}}{-(x+1-c)}).
	\end{align*}
	It is easy to verify that the first term is $ 0 $, and we now focus on the second term,
	\begin{align}
		\int_{a(c)}^{b(c)} f_{U}(x) d \tan ^{-1}(\frac{\sqrt{4 c-(x-1-c)^{2}}}{-(x+1-c)})
		=\int_{a(c)}^{b(c)} \frac{\log\left(1+x \right) }{1+\frac{4c-(x-1-c)^{2}}{(x+1-c)^{2}}}\cdot\frac{\sqrt{4c-(x-1-c)^{2}}+\frac{(x-1-c)(x+1-c)}{\sqrt{4c-(x-1-c)^{2}}}}{(x+1-c)^{2}}dx.\label{27}
	\end{align}
	By substituting $ x=1+c-2\sqrt{c}\cos(\theta) $, we obtain
	\begin{align}
		(\ref{27})=&\nonumber\frac{1}{2}\int_{0}^{2\pi}\left(\log\left(2+c-2\sqrt{c}\cos(\theta) \right)  \right)\frac{c-\sqrt{c}\cos(\theta) }{1+c-2\sqrt{c}\cos(\theta)} d\theta\\
		=&\nonumber\frac{1}{2}\int_{0}^{2\pi}[\log(1+\varrho\left(c \right)  )+\log(1+\sqrt{\tilde{c}}-2 \sqrt{\tilde{c}}\cos(\theta ))    ] \frac{c-\sqrt{c}\cos(\theta) }{1+c-2\sqrt{c}\cos(\theta)} d\theta\\
		=&\nonumber\frac{1}{2}\int_{0}^{2\pi}\log\left(1+\varrho\left(c \right) \right) \frac{c-\sqrt{c}\cos(\theta) }{1+c-2\sqrt{c}\cos(\theta)} d\theta+\nonumber\frac{1}{2}\int_{0}^{2\pi}\log(1+\sqrt{\tilde{c}}-2 \sqrt{\tilde{c}}\cos(\theta) ) \frac{c-\sqrt{c}\cos(\theta) }{1+c-2\sqrt{c}\cos(\theta)} d\theta.
	\end{align}
	For the first integral, by substituting $ \cos\theta=\frac{z+z^{-1}}{2} $, we turn it into a contour integral on $ \left|z \right|=1,  $
	\begin{align*}
		&\frac{1}{2}\int_{0}^{2\pi}\log\left(1+\varrho\left(c \right) \right) \frac{c-\sqrt{c}\cos(\theta) }{1+c-2\sqrt{c}\cos(\theta)} d\theta
		=\dfrac{\log\left(1+\varrho\left(c \right) \right) }{2}\oint_{\left| z\right|=1 }\dfrac{c-\sqrt{c}\frac{z+\frac{1}{z}}{2}   }{1+c-\sqrt{c}\left(z+\frac{1}{z} \right) }\frac{1}{iz}dz\\
		=&\dfrac{\log\left(1+\varrho\left(c \right) \right) }{2i}\oint_{\left| z\right|=1 }\dfrac{2cz-\sqrt{c}\left(z^{2}+1 \right) }{2z\left(z-\sqrt{c} \right)\left(1-\sqrt{c}z \right)  }dz.
	\end{align*}
	When $ c<1, $	$ 0 $ and $ \sqrt{c} $ are poles. The residues are $ \frac{1}{2} $ and $ -\frac{1}{2} $, respectively. By the residue theorem, the first integral is $ 0. $ 
	
	For the second integral, 	 
	\begin{align}
		&\nonumber\frac{1}{2}\int_{0}^{2\pi}\frac{\log\left(1+\tilde{c}-2\sqrt{\tilde{c}}\cos(\theta) \right)}{1+c-2\sqrt{c}\cos(\theta)}\left( c-\sqrt{c}\cos\theta\right)  d\theta\\
		=&\nonumber\frac{1}{2} \oint_{\left| z\right| =1} \log |1-\sqrt{\tilde{c}} z|^{2} \cdot \frac{c-\sqrt{c} \frac{z+z^{-1}}{2}}{1+c-2\sqrt{ c} \cdot \frac{z+z^{-1}}{2}} \frac{d z}{i z}
		=\frac{1}{4 i} \oint_{\left| z\right| =1} \log |1-\sqrt{\tilde{c}} z|^{2} \cdot \frac{2c z-\sqrt{c}\left(z^{2}+1\right)}{(z-\sqrt{c})(-\sqrt{c}z+1)z} d z\\
		=&\label{secint}\frac{1}{4 i} \oint_{\left| z\right| =1}\log\left( 1-\sqrt{\tilde{c}} z\right) \frac{2c z-\sqrt{c}\left(z^{2}+1\right)}{(z-\sqrt{c})(-\sqrt{c}z+1)z} d z+ \frac{1}{4 i} \oint_{\left| z\right| =1}\log\left( 1-\sqrt{\tilde{c}} \frac{1}{z}\right)\frac{2c z-\sqrt{c}\left(z^{2}+1\right)}{(z-\sqrt{c})(-\sqrt{c}z+1)z} dz.  
	\end{align}
	For the first term in (\ref{secint}),  when $ c<1 $, the pole is $\sqrt{c}$ , and the residue is $ -\log( 1-\sqrt{\tilde{c}c})  $. By using the residue theorem, the integral is $ -\frac{\pi}{2} \log( 1-\sqrt{\tilde{c}c}) $. The same argument also holds for the second term in (\ref{secint}), and the integral is also $ -\frac{\pi}{2} \log( 1-\sqrt{\tilde{c}c}) $ after some calculation. Therefore the second integral equals $ -\pi\log( 1-\sqrt{\tilde{c}c})  $. 
	Therefore $ \frac{M}{2\pi i}\oint_{\mathcal C}f_{U}\left(z \right)\frac{\underline{m}'(z)}{\underline{m}(z)}dz= M\log( 1-\sqrt{\tilde{c}c}), $
	 and the result is still valid if $ c $ is replaced by $ c_{nM} $.  Therefore, formula (\ref{Uextra}) holds.
	
\end{proof}

\begin{proof}[Proof of (\ref{Ui1}):] From Theorem A.1 in \cite{wang2013sphericity}, we have
	\begin{align*}
		I_{1}(f_{U} )=&\lim\limits_{r\downarrow1}\frac{1}{2\pi i} \oint_{\left| z\right|=1 }\log(1+\left|1+\sqrt{c}z \right|^{2}  )(\frac{z}{z^{2}-r^{-2}}-\frac{1}{z} )dz\\
		=&  \lim\limits_{r\downarrow1}\frac{1}{2\pi i} \oint_{\left| z\right|=1 }\log(1+\left|1+\sqrt{c}z \right|^{2}  )\frac{z}{z^{2}-r^{-2}}dz-\lim\limits_{r\downarrow1}\frac{1}{2\pi i} \oint_{\left| z\right|=1 }\log(1+\left|1+\sqrt{c}z \right|^{2}  )\frac{dz}{z}.
	\end{align*}
	For the first integral of $ I_{1}(f_{U} )  $, 
	\begin{align*}	
		&\lim\limits_{r\downarrow1}\frac{1}{2\pi i} \oint_{\left| z\right|=1 }\log(1+\left|1+\sqrt{c}z \right|^{2}  )\frac{z}{z^{2}-r^{-2}}dz
		=\lim\limits_{r\downarrow1}\frac{1}{2\pi i} \oint_{\left| z\right|=1 }\log\left[ (1+\varrho(c ) )(1+\sqrt{\tilde{c}}z )(1+\sqrt{\tilde{c}}\bar{z} ) \right] \frac{z}{z^{2}-r^{-2}}dz\\
		=&\lim\limits_{r\downarrow1}\frac{1}{2\pi i} \oint_{\left| z\right|=1 } \frac{\log (1+\varrho(c ) )z}{z^{2}-r^{-2}}dz+\lim\limits_{r\downarrow1}\frac{1}{2\pi i} \oint_{\left| z\right|=1 } \frac{\log (1+\sqrt{\tilde{c}}z )z}{z^{2}-r^{-2}}dz
		+\lim\limits_{r\downarrow1}\frac{1}{2\pi i} \oint_{\left| z\right|=1 }   \frac{\log(1+\sqrt{\tilde{c}}\frac{1}{z} )z}{z^{2}-r^{-2}}dz,
	\end{align*}
	where $ \tilde{c} $ and $ \varrho(c ) $ are defined in the proof of (\ref{Ucenter}). 
	For the first integral, the poles are $ \frac{1}{r} $ and $ -\frac{1}{r} $, the residues are both $ \frac{1}{2}\log(1+\varrho(c) )  $. Therefore, by the residue theorem, the integral is $ \log(1+\varrho(c) )  $. Similarly, for the second integral, the residues are $ \frac{1}{2}\log( 1+\frac{\sqrt{\tilde{c}}}{r})  $ and $ \frac{1}{2}\log( 1-\frac{\sqrt{\tilde{c}}}{r})  $. By the residue theorem, the integral is $ \frac{1}{2}\log(1-\tilde{c} ).  $ For the third integral,
	\begin{align*}
		&\lim\limits_{r\downarrow1}\frac{1}{2\pi i} \oint_{\left| z\right|=1 }\log(1+\sqrt{\tilde{c}}\frac{1}{z} )   \frac{z}{z^{2}-r^{-2}}dz\\
		=& \lim\limits_{r\downarrow1}\frac{1}{2\pi i} \oint_{\left| \xi\right|=1 }\log(1+\sqrt{\tilde{c}}\xi )   \frac{\frac{1}{\xi}}{\xi^{-2}-r^{-2}}\frac{1}{\xi^{2}}d\xi
		= \lim\limits_{r\downarrow1}r^{2}\frac{1}{2\pi i} \oint_{\left| \xi\right|=1 }\log( 1+\sqrt{\tilde{c}}\xi) \frac{1}{\xi(r+\xi )(r-\xi )  }d\xi,
	\end{align*}
	where the first integral results from the change of variable $ \xi=\frac{1}{z}. $ The poles are $ r $ and $ -r, $ and the residues are $ -\log(1+\sqrt{\tilde{c}}r )\frac{1}{2r^{2}}  $ and $ -\log(1-\sqrt{\tilde{c}}r)\frac{1}{2r^{2}} $, respectively. Then by the residue theorem, the integral is $ -\frac{1}{2}\log(1-\tilde{c} ).  $ 
	
	Collecting the three integral above leads to $$ \lim\limits_{r\downarrow1}\frac{1}{2\pi i} \oint_{\left| z\right|=1 }\log(1+\left|1+\sqrt{c}z \right|^{2}  )\frac{z}{z^{2}-r^{-2}}dz=\log(1+\varrho(c ) ).  $$

	For the second integral of $ I_{1}(f_{U} )  $, 
	\begin{align}
		&\nonumber\lim\limits_{r\downarrow1}\frac{1}{2\pi i} \oint_{\left| z\right|=1 }\log(1+\left|1+\sqrt{c}z \right|^{2}  )\frac{1}{z}dz
		=\lim\limits_{r\downarrow1}\frac{1}{2\pi i} \oint_{\left| z\right|=1 }\frac{\log(2+c )+ \log(1+\frac{\sqrt{c}}{2+c}(\frac{1}{z}+z )  )}{z}dz\\
		=&\lim\limits_{r\downarrow1}\frac{1}{2\pi i} \oint_{\left| z\right|=1 }\frac{\log(2+c )}{z}dz+\lim\limits_{r\downarrow1}\frac{1}{2\pi i} \oint_{\left| z\right|=1 }\frac{ \log(1+\frac{\sqrt{c}}{2+c}(\frac{1}{z}+z )  )}{z}dz.\label{i1int2}
	\end{align}	
	Since $ \left|\frac{\sqrt{c}}{2+c} (\frac{1}{z}+z)\right|= \frac{\sqrt{c}}{{2+c}}\left|\frac{1}{z}+z \right|   $, and $  \frac{\sqrt{c}}{{2+c}}< \frac{\sqrt{c}}{{1+c}}\leq \frac{1}{2} $, $ \left|\frac{1}{z}+z \right| \leq \left|\frac{1}{z} \right|+\left|z \right|=2,   $ then by using Taylor expansion, we have
	\begin{align*}
		& \oint_{\left| z\right|=1 }\frac{ \log(1+\frac{\sqrt{c}}{2+c}(\frac{1}{z}+z )  )}{z}dz\\
		=&\oint_{\left| z\right|=1 }\sum_{k=1}^{\infty}(-1)^{k+1}(\frac{\sqrt{c}}{2+c} )^{k}\frac{1}{k}(\frac{1}{z}+z)^{k}\frac{dz}{z} 
		=\sum_{k=1}^{\infty}(-1)^{k+1}(\frac{\sqrt{c}}{2+c} )^{k}\frac{1}{k}\oint_{\left| z\right|=1 }(\frac{1}{z}+z )^{k}\frac{dz}{z}\\
		=&\sum_{k=1}^{\infty}(-1)^{2k+1}(\frac{\sqrt{c}}{2+c} )^{2k}\frac{1}{2k}\oint_{\left| z\right|=1 }(\frac{1}{z}+z )^{2k}\frac{dz}{z} 
		=\sum_{k=1}^{\infty}(-1)(\frac{\sqrt{c}}{2+c} )^{2k}\frac{1}{2k}C_{2k}^{k}2\pi i
		= -\sum_{k=1}^{\infty}(\frac{\sqrt{c}}{2+c} )^{2k}\frac{(2k-1)!}{k!k!}2\pi i.
	\end{align*}
	Therefore the second integral equals $ \log(2+c )-\sum_{k=1}^{\infty}(\frac{\sqrt{c}}{2+c} )^{2k}\frac{(2k-1)!}{k!k!}.$ Thus $ I_{1}(f_{U} )=\log(1+\varrho )-  \log(2+c )+\sum_{k=1}^{\infty}(\frac{\sqrt{c}}{2+c} )^{2k}\frac{(2k-1)!}{k!k!}. $ Note that all the $ c $ in the formula above should be replaced by $ c_{nM} $ since the calculation is on the bulk part of LSS. The proof is finished.
\end{proof}
\begin{proof}[Proof of (\ref{Ui2}):]
	From Theorem A.1 in \cite{wang2013sphericity},  by Taylor expansion, we have $I_{2}(f_{U} )$ equals
	\begin{align*} 
		& \lim\limits_{r\downarrow1}\frac{1}{2\pi i} \oint_{\left| z\right|=1 }\log(1+\left|1+\sqrt{c}z \right|^{2}  )\frac{1}{z^{3}}dz 
		=\frac{1}{2\pi i} \oint_{\left| z\right|=1 }( \log(2+c) + \log(1+\frac{\sqrt{c}}{2+c}(\frac{1}{z}+z )  ) )\frac{1}{z^{3}}dz\\  
		&=\frac{1}{2\pi i} \oint_{\left| z\right|=1 }\log(1+\frac{\sqrt{c}}{2+c}(\frac{1}{z}+z )  )\frac{1}{z^{3}}dz
		=\frac{1}{2\pi i} \oint_{\left| z\right|=1 }\sum_{k=1}^{\infty}(-1)^{2k+1}(\frac{\sqrt{c}}{2+c} )^{2k}\frac{1}{2k}(\frac{1}{z}+z )^{2k}\frac{1}{z^{3}} dz\\
		&=\sum_{k=1}^{\infty}(-1)^{2k+1}(\frac{\sqrt{c}}{2+c} )^{2k}\frac{1}{2k}\frac{1}{2\pi i}\oint_{\left| z\right|=1 }(\frac{1}{z}+z )^{2k}\frac{1}{z^{3}} dz
		=\sum_{k=1}^{\infty}(-1)^{2k+1}(\frac{\sqrt{c}}{2+c} )^{2k}\frac{1}{2k}C_{2k}^{k-1}
		=-\sum_{k=1}^{\infty}(\frac{\sqrt{c}}{2+c} )^{2k}\frac{(2k-1)!}{(k-1)!(k+1)!}.
	\end{align*}
	Notice that all the $ c $ in the formula above should be replaced by $ c_{nM} $ since the calculation is on the bulk part of LSS. The proof is finished.
\end{proof}
\begin{proof}[Proof of (\ref{Uj1}):]  From Theorem A.1 in \cite{wang2013sphericity}, we have $J_{1}( f_{U},f_{U})$ equals
	\begin{align*}
		&\lim\limits_{r\downarrow1}-\frac{1}{4 \pi^{2}} \oint_{\left| z_{1}\right|=1 }\oint_{\left| z_{2}\right|=1 }\frac{\log(1+\left|1+\sqrt{c}z_{1} \right|^{2}  ) \log(1+\left|1+\sqrt{c}z_{2} \right|^{2}  )}{( z_{1}-rz_{2})^{2} }dz_{1}dz_{2}\\
		&=\lim\limits_{r\downarrow1}-\frac{1}{4 \pi^{2}}\oint_{\left| z_{2}\right|=1 } \log(1+\left|1+\sqrt{c}z_{2} \right|^{2}  )\oint_{\left| z_{1}\right|=1 }\frac{\log(1+\left|1+\sqrt{c}z_{1} \right|^{2}  )}{( z_{1}-rz_{2})^{2}}dz_{1}dz_{2}.
	\end{align*}
	Since $ r>1, $ thus $ rz_{2} $ is not a pole.
	\begin{align*}
		&\oint_{\left| z_{1}\right|=1 }\frac{\log(1+\left|1+\sqrt{c}z_{1} \right|^{2}  )}{( z_{1}-rz_{2})^{2}}dz_{1}
		=\oint_{\left| z_{1}\right|=1 }\frac{\sum_{k=1}^{\infty}(-1)^{k+1}(\frac{\sqrt{c}}{2+c} )^{k}\frac{1}{k}(\frac{1}{z_{1}}+z_{1} )^{k}}{(z_{1}-rz_{2} )^{2} }dz_{1}\\
		&=\oint_{\left| z_{1}\right|=1 }\frac{\sum_{k=1}^{\infty}(-1)^{2k}(\frac{\sqrt{c}}{2+c} )^{2k-1}\frac{1}{2k-1}(\frac{1}{z_{1}}+z_{1} )^{2k-1}}{(z_{1}-rz_{2} )^{2} }dz_{1}\\
		&=\sum_{k=1}^{\infty}(\frac{\sqrt{c}}{2+c} )^{2k-1}\frac{1}{2k-1}\oint_{\left| z_{1}\right|=1 }(\frac{1}{z_{1}}+z_{1} )^{2k-1}\frac{1}{(z_{1}-rz_{2} )^{2}}dz_{1}\\
		&=\sum_{k=1}^{\infty}(\frac{\sqrt{c}}{2+c} )^{2k-1}\frac{1}{2k-1}2\pi i\frac{(2k-1)!}{k!(k-1)!}\frac{1}{r^{2}z_{2}^{2}}
		=\sum_{k=1}^{\infty}(\frac{\sqrt{c}}{2+c} )^{2k-1}2\pi i\frac{(2k-2)!}{k!(k-1)!}\frac{1}{r^{2}z_{2}^{2}}. 
	\end{align*}	
	Then we have $J_{1}( f_{U}, f_{U})$ equals
	\begin{align*}
		\lim\limits_{r\downarrow1}-\frac{1}{4 \pi^{2}}\sum_{k=1}^{\infty}(\frac{\sqrt{c}}{2+c} )^{2k-1}2\pi i\frac{(2k-2)!}{k!(k-1)!}\oint_{\left| z_{2}\right|=1 } \log(1+\left|1+\sqrt{c}z_{2} \right|^{2}  )\frac{1}{r^{2}z_{2}^{2}}dz_{2},
	\end{align*}
	and by using Taylor expansion to $\log(1+\left|1+\sqrt{c}z_{2} \right|^{2} ) $, we can obtain that 
	$
	J_{1}( f_{U}, f_{U})	=( \sum_{k=1}^{\infty}(\frac{\sqrt{c}}{2+c} )^{2k-1}\frac{(2k-2)!}{k!(k-1)!})^{2} 
	$	
	(the contour integral about $ z_{2} $ is handled the same way as $ z_{1} $). 
	Similarly, from Theorem A.1 in \cite{wang2013sphericity}, we have $J_{2}(f_{U}, f_{U} )$ equals
	\begin{align*}
		\lim\limits_{r\downarrow1}-\frac{1}{4 \pi^{2}} \oint_{\left| z_{1}\right|=1 }\frac{ \log(1+\left|1+\sqrt{c}z_{1} \right|^{2}  ) }{z_{1}^{2}}dz_{1}\oint_{\left| z_{2}\right|=1 }\frac{ \log(1+\left|1+\sqrt{c}z_{2} \right|^{2}  ) }{z_{2}^{2}}dz_{2}.
	\end{align*}	
	By Taylor expansion, we obtain 
	\begin{align*}
		&\oint_{\left| z_{1}\right|=1 }\frac{ \log(1+\left|1+\sqrt{c}z_{1} \right|^{2}  ) }{z_{1}^{2}}dz_{1}
		=\oint_{\left| z_{1}\right|=1 }\frac{\log(2+c) }{z_{1}^{2}}dz_{1}+
		\oint_{\left| z_{1}\right|=1 }\frac{ \log(1+\frac{\sqrt{c}}{2+c}(\frac{1}{z_{1}}+z_{1} )  )}{z_{1}^{2}}dz_{1}\\
		&=\sum_{k=1}^{\infty}(-1 )^{k+1}(\frac{\sqrt{c}}{2+c} )^{k}\frac{1}{k}\oint_{\left| z_{1}\right|=1 }(\frac{1}{z_{1}}+z_{1} )^{k}\frac{1}{z_{1}^{2}}dz_{1} 
		= \sum_{k=1}^{\infty}(-1 )^{2k}(\frac{\sqrt{c}}{2+c} )^{2k-1}\frac{1}{2k-1}\oint_{\left| z_{1}\right|=1 }(\frac{1}{z_{1}}+z_{1} )^{2k-1}\frac{1}{z_{1}^{2}}dz_{1} \\&=2\pi i \sum_{k=1}^{\infty}(\frac{\sqrt{c}}{2+c} )^{2k-1}\frac{(2k-2 )! }{k!(k-1 )! },
	\end{align*}	
	thus 
$
		J_{2}(f_{U}, f_{U} )=J_{1}(f_{U}, f_{U} )=(\sum_{k=1}^{\infty}(\frac{\sqrt{c}}{2+c} )^{2k-1}\frac{(2k-2)!}{k!( k-1) !})^{2}.
$	Notice that all the $ c $ in the formula above should be replaced by $ c_{nM} $ since the calculation is on the bulk part of LSS. The proof is finished.
\end{proof}
\begin{proof}[Proof of (\ref{Wextra}):] Similarly to proof of (\ref{Uextra}), we have 
	\begin{align*} &\oint_{\mathcal C}f_{W}(z )\frac{\underline{m}'(z)}{\underline{m}(z)}dz
		=2i 	\int_{a(c)}^{b(c)}f_{W}^{'}(x ) \Im \log\underline{m}(x+i\varepsilon)dx
		=-2i\int_{a(c)}^{b(c)}xd \tan ^{-1}(\frac{\sqrt{4 c-(x-1-c)^{2}}}{-(x+1-c)}).
	\end{align*} 
	Since
	\begin{align*}
		&\int_{a(c)}^{b(c)}xd \tan ^{-1}(\frac{\sqrt{4 c-(x-1-c)^{2}}}{-(x+1-c)})
		=\int_{a(c)}^{b(c)} \frac{x }{1+\frac{4c-(x-1-c)^{2}}{(x+1-c)^{2}}}\cdot\frac{\sqrt{4c-(x-1-c)^{2}}+\frac{(x-1-c)(x+1-c)}{\sqrt{4c-(x-1-c)^{2}}}}{(x+1-c)^{2}}dx.
	\end{align*}	
	By substituting $ x=1+c-2\sqrt{c}\cos(\theta) $, it equals
	\begin{align*}
		&\frac{1}{2}\int_{0}^{2\pi}(1+c-2\sqrt{c}\cos(\theta) )  \frac{c-\sqrt{c}\cos(\theta) }{1+c-2\sqrt{c}\cos(\theta)} d\theta
		=\int_{0}^{\pi}( c-\sqrt{c}\cos(\theta))  d\theta=\pi c.
	\end{align*}	
	Therefore $  \oint_{\mathcal C}f_{W}(z )\frac{\underline{m}'(z)}{\underline{m}(z)}dz =-2\pi i c,$ then $ \frac{M}{2\pi i}\oint_{\mathcal C}f_{W}(z )\frac{\underline{m}'(z)}{\underline{m}(z)}dz=-Mc, $
	the proof is finished.
\end{proof}

\begin{proof}[Proof of (\ref{Vcenter}):] Since $ f_{V}( x)=\frac{x}{1+x},  $ then
	\begin{align*}
		\int f_{V}( x)dF^{c_{n},H_{n}}(x)=\int_{a(c_{n})}^{b(c_{n})}\frac{x}{1+x}\frac{1}{2\pi xc_{n}}\sqrt{( b(c_{n})-x)( x-a(c_{n}))  }dx.
	\end{align*}	
	By using the variable change $ x=1+c_{n}-2\sqrt{c_{n}}\cos(\theta), $ $ 0\leq\theta\leq\pi, $ we have $\int f_{V}( x)dF^{c_{n},H_{n}}(x)$ equals
	\begin{align*}
		&\frac{1}{2\pi c_{n}}\int_{0}^{\pi}\frac{1}{2+c_{n}-2\sqrt{c_{n}}\cos(\theta)}4c_{n}\sin^{2}(\theta)d\theta
		=\frac{1}{4\pi c_{n}}\int_{0}^{2\pi}\frac{1}{(1+\varrho(c_{n} ) ) ( 1+\tilde{c}_{n}-2\sqrt{\tilde{c}_{n}} \cos( \theta) ) }4c_{n}\sin^{2}(\theta)d\theta \\
		=&\frac{1}{4\pi c_{n}}\frac{1}{1+\varrho(c_{n} ) }\int_{0}^{2\pi}\frac{1}{ 1+\tilde{c}_{n}-2\sqrt{\tilde{c}_{n}} \cos( \theta)}4c_{n}\sin^{2}(\theta)d\theta 
		=-\frac{1}{4\pi i}\frac{1}{1+\varrho(c_{n} ) }\oint_{\left| z\right|=1 }\frac{(z-1)^{2}(z+1)^{2}}{z^{2}( z-\sqrt{\tilde{c}_{n}})(1-\sqrt{\tilde{c}_{n}}z )  }dz.
	\end{align*}	
	When $ \tilde{c}_{n}<1, $ the poles are $ 0 $ and $ \sqrt{\tilde{c}_{n}}. $ The residues are $ -\frac{1+\tilde{c}_{n}}{\tilde{c}_{n}} $ and $ \frac{1-\tilde{c}_{n}}{\tilde{c}_{n}}, $ respectively. Then by the residue theorem, $ -\frac{1}{4\pi i}\frac{1}{1+\varrho(c_{n} ) }\oint_{\left| z\right|=1 }\frac{(z-1)^{2}(z+1)^{2}}{z^{2}( z-\sqrt{\tilde{c}_{n}})(1-\sqrt{\tilde{c}_{n}}z )  }dz $ equals $ \frac{1}{1+\varrho(c_{n} ) }. $ Notice that all the $ c_{n} $ in the formulas above should be replaced by $ c_{nM} $. The proof is finished. 
\end{proof}
\begin{proof}[Proof of (\ref{Vextra}):] Similarly to proof of (\ref{Uextra}), we have 
	\begin{align*} \oint_{\mathcal C}f_{V}(z )\frac{\underline{m}'(z)}{\underline{m}(z)}dz
		=2i 	\int_{a(c)}^{b(c)}f_{V}^{'}(x ) \Im \log\underline{m}(x+i\varepsilon)dx
		=-2i\int_{a(c)}^{b(c)}\frac{x}{1+x}d \tan ^{-1}(\frac{\sqrt{4 c-(x-1-c)^{2}}}{-(x+1-c)}).
	\end{align*} 
	Since
	\begin{align*}
		\int_{a(c)}^{b(c)}\frac{x}{1+x}d \tan ^{-1}(\frac{\sqrt{4 c-(x-1-c)^{2}}}{-(x+1-c)})
		=\frac{1}{2}\int_{0}^{2\pi}\frac{(c-\sqrt{c}\cos(\theta ))d\theta }{2+c-2\sqrt{c}\cos(\theta ) }
		=\int_{0}^{\pi}\frac{d\theta}{2}+ ( \frac{c}{2}-1) \int_{0}^{\pi}\frac{d\theta}{2+c-2\sqrt{c}\cos(\theta )}.
	\end{align*}	
	The first integral is $ \frac{\pi}{2}.  $ For the second integral, it equals
	\begin{align*}
		( \frac{c}{2}-1)\frac{1}{1+\varrho}\frac{1}{2}\int_{0}^{\pi}\frac{1}{ 1+\tilde{c}-2\sqrt{\tilde{c}}\cos(\theta) }d\theta
		= \frac{c-2}{4i}\frac{1}{1+\varrho }\oint_{\left| z\right|=1 }\frac{1}{( z-\sqrt{\tilde{c}}) ( 1-\sqrt{\tilde{c}}z) }dz.
	\end{align*}	 
	When $ \tilde{c}<1, $ the pole is $ \sqrt{\tilde{c}}, $ and the residue is $ \frac{1}{1-\tilde{c}}. $ By using the residue theorem, the integral \\$ \frac{c-2}{4i(1+\varrho)}\oint_{\left| z\right|=1 }\frac{1}{( z-\sqrt{\tilde{c}}) ( 1-\sqrt{\tilde{c}}z) }dz=\frac{\pi(c-2 ) }{2( 1+\varrho)(1-\tilde{c} )  }, $ then we have $ \int_{a(c)}^{b(c)}\frac{x}{1+x}d\tan ^{-1}(\frac{\sqrt{4 c-(x-1-c)^{2}}}{-(x+1-c)}) = \frac{\pi}{2}+\frac{\pi(c-2 ) }{2( 1+\varrho)(1-\tilde{c} )  }, $ therefore $ \oint_{\mathcal C}f_{V}(z )\frac{\underline{m}'(z)}{\underline{m}(z)}dz $ equals $ -\frac{1}{2}-\frac{c-2  }{2( 1+\varrho)(1-\tilde{c} )  }, $ then the proof is finished.
\end{proof}
\begin{proof}[Proof of (\ref{Vi1}):] From Theorem A.1 in \cite{wang2013sphericity}, we have
	\begin{align*}
		&I_{1}(f_{V} )=\lim\limits_{r\downarrow1}\frac{1}{2\pi i} \oint_{\left| z\right|=1 }\frac{\left|1+\sqrt{c}z \right|^{2}}{1+\left|1+\sqrt{c}z \right|^{2} } (\frac{z}{z^{2}-r^{-2}}-\frac{1}{z} )dz\\
		&=  \lim\limits_{r\downarrow1}\frac{1}{2\pi i} \oint_{\left| z\right|=1 }\frac{\left|1+\sqrt{c}z \right|^{2} }{1+\left|1+\sqrt{c}z \right|^{2} } \frac{z}{z^{2}-r^{-2}}dz-\lim\limits_{r\downarrow1}\frac{1}{2\pi i} \oint_{\left| z\right|=1 }\frac{\left|1+\sqrt{c}z \right|^{2} }{1+\left|1+\sqrt{c}z \right|^{2} }\frac{1}{z}dz. 
	\end{align*}
	For the first integral, 
	\begin{align}
		\lim\limits_{r\downarrow1}\frac{1}{2\pi i} \oint_{\left| z\right|=1 }\frac{\left|1+\sqrt{c}z \right|^{2} }{1+\left|1+\sqrt{c}z \right|^{2} } \frac{z}{z^{2}-r^{-2}}dz
		=\lim\limits_{r\downarrow1}\frac{1}{2\pi i} \oint_{\left| z\right|=1 } \frac{z}{z^{2}-r^{-2}}dz-\lim\limits_{r\downarrow1}\frac{1}{2\pi i} \oint_{\left| z\right|=1 }\frac{1}{1+\left|1+\sqrt{c}z \right|^{2} } \frac{z}{z^{2}-r^{-2}}dz\label{first term of Vi1}.
	\end{align}	
	For the first term of (\ref{first term of Vi1}), by using the residue theorem, the integral is $ 1. $  For the second integral,  
	\begin{align}
		&\nonumber\lim\limits_{r\downarrow1}\frac{1}{2\pi i} \oint_{\left| z\right|=1 }\frac{1}{1+\left|1+\sqrt{c}z \right|^{2} } \frac{z}{z^{2}-r^{-2}}dz
		=\lim\limits_{r\downarrow1}\frac{1}{2\pi i} \oint_{\left| z\right|=1 }\frac{1}{(1+\varrho )( 1+\sqrt{\tilde{c}}z)(1+\sqrt{\tilde{c}}\frac{1}{z} )      }\cdot\frac{z}{(z-\frac{1}{r} )(z+\frac{1}{r} )  }dz\\
		=& \lim\limits_{r\downarrow1}\frac{1}{2\pi i}\frac{1}{1+\varrho} \oint_{\left| z\right|=1 }\frac{z^{2}}{ (1+\sqrt{\tilde{c}}z     )( z+\sqrt{\tilde{c}} )( z+\frac{1}{r})(z-\frac{1}{r} )    }dz. \label{666}
	\end{align}	
	For $ \oint_{\left| z\right|=1 }\frac{z^{2}}{ (1+\sqrt{\tilde{c}}z     )( z+\sqrt{\tilde{c}} )( z+\frac{1}{r})(z-\frac{1}{r} )    }dz $,  it has $ -\sqrt{\tilde{c}}$, $ -\frac{1}{r}, $ $ \frac{1}{r} $ three poles, the residues are $ \frac{\tilde{c}}{(1-\tilde{c} )(-\sqrt{\tilde{c}}+ \frac{1}{r})  (-\sqrt{\tilde{c}}- \frac{1}{r}) }, $ $ \frac{1/r^{2}}{ (1-\frac{\sqrt{\tilde{c}}}{r} )( \sqrt{\tilde{c}}-\frac{1}{r})\cdot\frac{-2}{r}   } $, $ \frac{1/r^{2}}{ (1+\frac{\sqrt{\tilde{c}}}{r} ) ( \sqrt{\tilde{c}}+\frac{1}{r})\cdot\frac{2}{r} }. $ Then the summation of residues tend to $ \frac{\tilde{c}}{-( \tilde{c}-1)^{2} }+\frac{1}{2(\sqrt{\tilde{c}}-1 )^{2} } +\frac{1}{2(\sqrt{\tilde{c}}+1 )^{2} }.$ Therefore the integral (\ref{666}) equals $ \frac{1}{1+\varrho} ( \frac{\tilde{c}}{-( \tilde{c}-1)^{2} }+\frac{1}{2(\sqrt{\tilde{c}}-1 )^{2} } +\frac{1}{2(\sqrt{\tilde{c}}+1 )^{2} }). $ Then the equation $ (\ref{first term of Vi1}) $ equals $ 1- \frac{1}{1+\varrho} ( \frac{\tilde{c}}{-( \tilde{c}-1)^{2} }+\frac{1}{2(\sqrt{\tilde{c}}-1 )^{2} } +\frac{1}{2(\sqrt{\tilde{c}}+1 )^{2} }). $ 
	
	Then we consider the second integral of $ I_{1}( f_{V}).  $
	\begin{align*}
		&\lim\limits_{r\downarrow1}\frac{1}{2\pi i} \oint_{\left| z\right|=1 }\frac{\left|1+\sqrt{c}z \right|^{2} }{1+\left|1+\sqrt{c}z \right|^{2} }\frac{1}{z}dz
		=\lim\limits_{r\downarrow1}\frac{1}{2\pi i} \oint_{\left| z\right|=1 }\frac{1}{z}dz- \lim\limits_{r\downarrow1}\frac{1}{2\pi i} \oint_{\left| z\right|=1 }\frac{1 }{1+\left|1+\sqrt{c}z \right|^{2} }\frac{1}{z}dz.
	\end{align*}	
	The first integral is $ 1. $ By Taylor expansion, the second integral equals
	\begin{align*}
		&\lim\limits_{r\downarrow1}\frac{1}{2\pi i}\frac{1}{2+c}\sum_{k=0}^{\infty}( -1)^{2k}( \frac{\sqrt{c}}{2+c})^{2k}  \oint_{\left| z\right|=1 }(\frac{1}{z}+z )^{2k}\frac{1}{z}dz
		=\frac{1}{2+c}\sum_{k=0}^{\infty}( \frac{\sqrt{c}}{2+c})^{2k} \frac{( 2k) !}{k!k!},
	\end{align*}	
	therefore the second integral of $ I_{1}( f_{V})  $ equals $ 1-\frac{1}{2+c}\sum_{k=0}^{\infty}( \frac{\sqrt{c}}{2+c})^{2k} \frac{( 2k) !}{k!k!}. $ Collecting all the integrals of $  I_{1}( f_{V}) $, it equals $\frac{1}{2+c}\sum_{k=0}^{\infty}( \frac{\sqrt{c}}{2+c})^{2k} \frac{( 2k) !}{k!k!}-  \frac{1}{1+\varrho}\cdot( \frac{\tilde{c}}{-( \tilde{c}-1)^{2} }+\frac{1}{2(\sqrt{\tilde{c}}-1 )^{2} } +\frac{1}{2(\sqrt{\tilde{c}}+1)^{2} }).$
	Notice that all the $ c $ in the formulas above should be replaced by $ c_{nM} $ since the calculation is on the bulk part of LSS. The proof is finished.
\end{proof}
\begin{proof}[Proof of (\ref{Vi2}):] From Theorem A.1 in \cite{wang2013sphericity}, we have $I_{2}(f_{V} )$ equals
	\begin{align*}
		\frac{1}{2\pi i} \oint_{\left| z\right|=1 }\frac{\left|1+\sqrt{c}z \right|^{2}}{1+\left|1+\sqrt{c}z \right|^{2} } \frac{1}{z^{3}} dz
		=\frac{1}{2\pi i} \oint_{\left| z\right|=1 } \frac{1}{z^{3}} dz-\frac{1}{2\pi i} \oint_{\left| z\right|=1 }\frac{1}{1+\left|1+\sqrt{c}z \right|^{2} } \frac{1}{z^{3}} dz
		=-\frac{1}{2\pi i} \oint_{\left| z\right|=1 }\frac{1}{1+\left|1+\sqrt{c}z \right|^{2} } \frac{1}{z^{3}} dz.
	\end{align*}
	By Taylor expansion, $\frac{1}{1+\left|1+\sqrt{c}z \right|^{2} } =\frac{1}{2+c}\sum_{k=0}^{\infty}(-1 )^{k} (\frac{\sqrt{c}}{2+c})^{k}(\frac{1}{z}+z )^{k},   $ then $I_{2}(f_{V})  $ equals
	\begin{align*}
		&	-\frac{1}{2\pi i} \oint_{\left| z\right|=1 }\frac{1}{2+c}\sum_{k=0}^{\infty}(-1 )^{k} (\frac{\sqrt{c}}{2+c} )^{k}(\frac{1}{z}+z )^{k}\frac{1}{z^{3}} dz
		=	-\frac{1}{2\pi i} \frac{1}{2+c}\sum_{k=0}^{\infty}(-1 )^{k} (\frac{\sqrt{c}}{2+c} )^{k}\oint_{\left| z\right|=1 }(\frac{1}{z}+z )^{k}\frac{1}{z^{3}} dz\\
		&=-\frac{1}{2\pi i} \frac{1}{2+c}\sum_{k=0}^{\infty}(-1 )^{2k} (\frac{\sqrt{c}}{2+c} )^{2k}\oint_{\left| z\right|=1 }(\frac{1}{z}+z )^{2k}\frac{1}{z^{3}} dz
		=-\frac{1}{2+c}\sum_{k=0}^{\infty}(\frac{\sqrt{c}}{2+c} )^{2k}\frac{(2k)!}{(k-1)!(k+1)!}
	\end{align*}	
	Notice that all the $ c $ in the formulas above should be replaced by $ c_{nM} $ since the calculation is on the bulk part of LSS. The proof is finished.
\end{proof}
\begin{proof}[Proof of (\ref{Vj1}):] From Theorem A.1 in \cite{wang2013sphericity}, we have $J_{1}(f_{V},f_{V} )$ equals
	\begin{align*}
		 \lim\limits_{r\downarrow1}-\frac{1}{4\pi^{2}}\oint_{\left| z_{1}\right|=1 }\oint_{\left| z_{2}\right|=1 }\frac{\frac{\left|1+\sqrt{c}z_{1} \right|^{2}}{1+\left|1+\sqrt{c}z_{1} \right|^{2} }\frac{\left|1+\sqrt{c}z_{2} \right|^{2}}{1+\left|1+\sqrt{c}z_{2} \right|^{2} }}{(z_{1}-rz_{2} )^{2} }dz_{1}dz_{2}
		= \lim\limits_{r\downarrow1}-\frac{1}{4\pi^{2}}\oint_{\left| z_{2}\right|=1 }\frac{\left|1+\sqrt{c}z_{2} \right|^{2}}{1+\left|1+\sqrt{c}z_{2} \right|^{2} }\oint_{\left| z_{1}\right|=1 }\frac{ \frac{\left|1+\sqrt{c}z_{1} \right|^{2}}{1+\left|1+\sqrt{c}z_{1} \right|^{2} }}{(z_{1}-rz_{2} )^{2}}dz_{1}dz_{2}.
	\end{align*}	
	For the integral $ \oint_{\left| z_{1}\right|=1 }\frac{ \frac{\left|1+\sqrt{c}z_{1} \right|^{2}}{1+\left|1+\sqrt{c}z_{1} \right|^{2} }}{(z_{1}-rz_{2} )^{2}}dz_{1}, $  it equals
	\begin{align*}
		&\oint_{\left| z_{1}\right|=1 }\frac{ 1}{ ( z_{1}-rz_{2})^{2}  }dz_{1}-\oint_{\left| z_{1}\right|=1 }\frac{1}{1+\left|1+\sqrt{c}z_{2} \right|^{2} }\frac{ 1}{ ( z_{1}-rz_{2})^{2}  }dz_{1}\\
		=&-\oint_{\left| z_{1}\right|=1 }\frac{1}{1+\left|1+\sqrt{c}z_{2} \right|^{2} }\frac{ 1}{ ( z_{1}-rz_{2})^{2}  }dz_{1}
		= -\oint_{\left| z_{1}\right|=1 }\frac{1}{2+c}\sum_{k=0}^{\infty}(-1 )^{k} (\frac{\sqrt{c}}{2+c} )^{k}(\frac{1}{z_{1}}+z_{1} )^{k}\frac{ 1}{ ( z_{1}-rz_{2})^{2}  }dz_{1}\\
		=&-\frac{1}{2+c}\sum_{k=0}^{\infty}(-1 )^{k} (\frac{\sqrt{c}}{2+c} )^{k}\oint_{\left| z_{1}\right|=1 }(\frac{1}{z_{1}}+z_{1} )^{k}\frac{ 1}{ ( z_{1}-rz_{2})^{2}  }dz_{1}
		=\frac{1}{2+c}\sum_{k=0}^{\infty}(\frac{\sqrt{c}}{2+c} )^{2k+1}\frac{(2k+1)!}{(k+1)!k!}\frac{1}{r^{2}z_{2}^{2}}2\pi i.
	\end{align*}	
	By using the same methods as above, then $  J_{1}(f_{V},f_{V} )$ equals
	\begin{align*}
		&\lim\limits_{r\downarrow1}-\frac{1}{4\pi^{2} }\oint_{\left| z_{2}\right|=1 }\frac{\left|1+\sqrt{c}z_{2} \right|^{2}}{1+\left|1+\sqrt{c}z_{2} \right|^{2} }\frac{1}{2+c}\sum_{k=0}^{\infty}(\frac{\sqrt{c}}{2+c} )^{2k+1}\frac{(2k+1)!}{(k+1)!k!}\frac{1}{r^{2}z_{2}^{2}}2\pi idz_{2}\\
		=&\lim\limits_{r\downarrow1}-\frac{1}{4\pi^{2} }2\pi i\frac{1}{2+c}\sum_{k=0}^{\infty}(\frac{\sqrt{c}}{2+c} )^{2k+1}\frac{(2k+1)!}{(k+1)!k!}\oint_{\left| z_{2}\right|=1 }\frac{\left|1+\sqrt{c}z_{2} \right|^{2}}{1+\left|1+\sqrt{c}z_{2} \right|^{2} }\frac{1}{r^{2}z_{2}^{2}}dz_{2}\\
		=&( \frac{1}{2+c}\sum_{k=0}^{\infty}(\frac{\sqrt{c}}{2+c} )^{2k+1}\frac{(2k+1)!}{(k+1)!k!})^{2}. 
	\end{align*}	
	Then we consider $ J_{2}(f_{V},f_{V} ). $ Since
$
		J_{2}(f_{V},f_{V} )=-\frac{1}{4\pi^{2}}\oint_{\left| z_{1}\right|=1 }\frac{\frac{\left|1+\sqrt{c}z_{1} \right|^{2}}{1+\left|1+\sqrt{c}z_{1} \right|^{2}}}{z_{1}^{2}}dz_{1} \oint_{\left| z_{2}\right|=1 }\frac{\frac{\left|1+\sqrt{c}z_{2} \right|^{2}}{1+\left|1+\sqrt{c}z_{2} \right|^{2}}}{z_{2}^{2}}dz_{2}. 
$	
	For the integral \\$ \oint_{\left| z_{1}\right|=1 }\frac{\frac{\left|1+\sqrt{c}z_{1} \right|^{2}}{1+\left|1+\sqrt{c}z_{1} \right|^{2}}}{z_{1}^{2}}dz_{1}, $ by Taylor expansion, it equals
	\begin{align*}
		&-\frac{1}{2+c}\sum_{k=0}^{\infty}(-1)^{k}( \frac{\sqrt{c}}{2+c})^{k} \oint_{\left| z_{1}\right|=1 }( \frac{1}{z_{1}}+z_{1})^{k}\frac{1}{z_{1}^{2}} dz_{1}\\
		=&-\frac{1}{2+c}\sum_{k=0}^{\infty}(-1)^{2k+1}( \frac{\sqrt{c}}{2+c})^{2k+1} \oint_{\left| z_{1}\right|=1 }( \frac{1}{z_{1}}+z_{1})^{2k+1}\frac{1}{z_{1}^{2}} dz_{1}
		=\frac{1}{2+c}\sum_{k=0}^{\infty}( \frac{\sqrt{c}}{2+c})^{2k+1}\frac{(2k+1)!}{k!(k+1)!}2\pi i.
	\end{align*}	
	Therefore $ J_{2}(f_{V},f_{V} )=J_{1}(f_{V},f_{V} )=(\frac{1}{2+c}\sum_{k=0}^{\infty}( \frac{\sqrt{c}}{2+c})^{2k+1}\frac{(2k+1)!}{k!(k+1)!})^{2}.  $
	Notice that all the $ c $ in the formulas above should be replaced by $ c_{nM} $ since the calculation is on the bulk part of LSS. The proof is finished.
\end{proof}

\section*{Acknowledgments}

Zhijun Liu was partially supported by Fundamental Research Funds for Northeastern University No.N25X\\QD026. 
Jiang Hu was partially supported by NSFC Grants No. 12171078, No. 12292980, No. 12292982, National Key R $\&$ D Program of China No. 2020YFA0714102, and Fundamental Research Funds for the Central Universities No. 2412023YQ003. Zhidong Bai was partially supported by NSFC Grants No.12171198, No.12271536. Zhihui Lv was supported by Guangdong Basic and Applied Basic Research Foundation (2022A1515110084).

\section*{Supplementary material}
Supplementary material includes  simulation results.




\newpage
\begin{center}
	{\bf \Large	Supplementary material for ``asymptotic distributions of four linear hypotheses test statistics under generalized spiked model"}
\end{center}

\begin{center}
	Zhijun Liu, Jiang Hu, Zhidong Bai, Zhihui Lv
\end{center}

\

\

In this document we present some  simulation results  involved in \cite{Liu25}. The number of schemes (equations, theorems, lemmas, etc.) is shared with the main document so that there are no misunderstandings with the use of references.

\section{Simulation results}
\label{app1}
In this document we  present some comparisons between empirical distributions of $U,W,V$ and standard normal curves (red lines) under Models 1--4 when samples are from $Dt_1$(N(0,1)) and $Dt_2$ (Gamma (4,0.5)-2), respectively.

\begin{figure}[htbp]
	\centering
	\subfigure[Model 1]{
		\includegraphics[width=3.5cm,height=3.5cm]{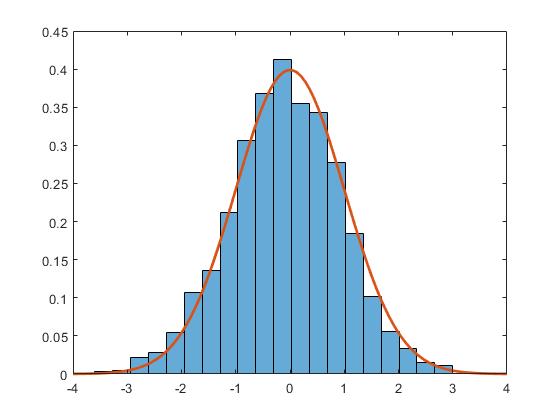}\label{u1}}
	\subfigure[Model 2]{
		\includegraphics[width=3.5cm,height=3.5cm]{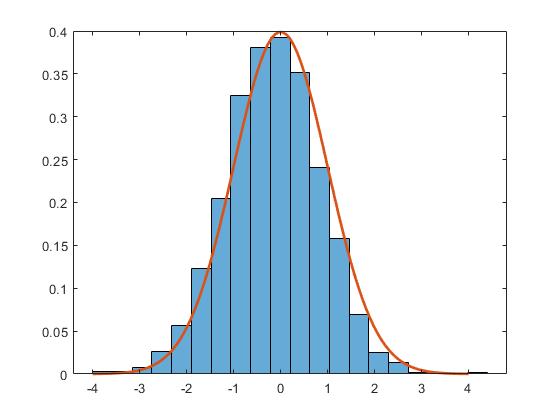}\label{u2}}
	\subfigure[Model 3]{
		\includegraphics[width=3.5cm,height=3.5cm]{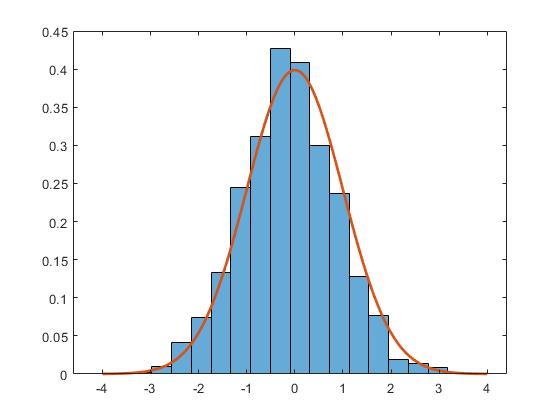}\label{u3}}
	\subfigure[Model 4]{ 
		\includegraphics[width=3.5cm,height=3.5cm]{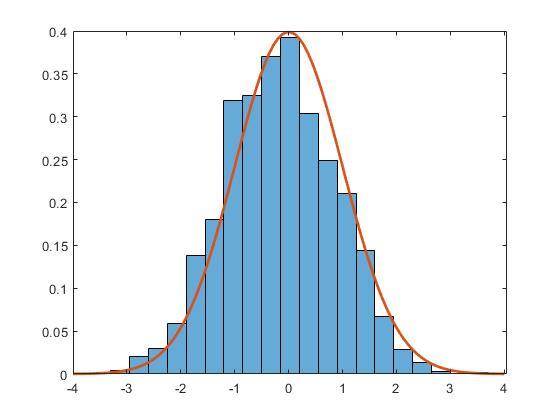}\label{u4}}
	
	\caption{
		Comparisons between empirical distributions of statistics $U$ and standard normal curves	under Models 1--4, respectively, when samples are from $Dt_1$	
	}
	\label{fig1}
\end{figure}

\begin{figure}[htbp]
	\centering
	\subfigure[Model 1]{
		\includegraphics[width=3.5cm,height=3.5cm]{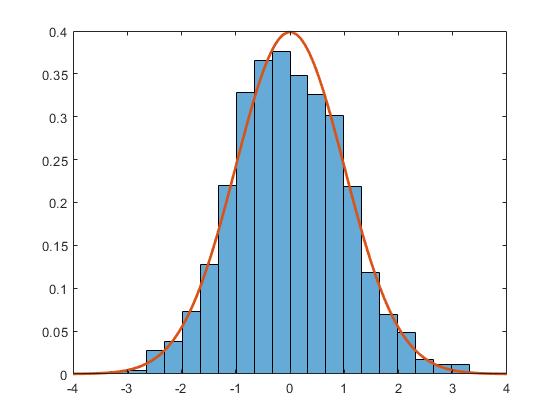}\label{w1}}
	\subfigure[Model 2]{
		\includegraphics[width=3.5cm,height=3.5cm]{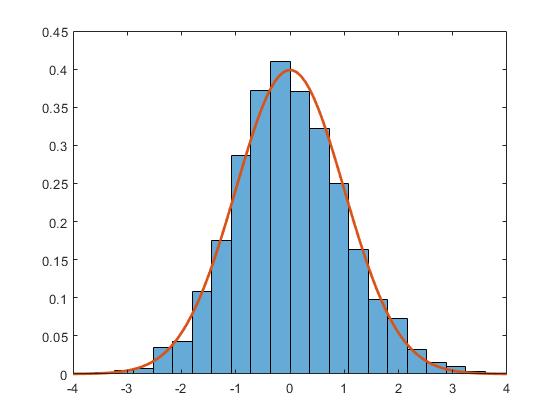}\label{w2}}
	\subfigure[Model 3]{
		\includegraphics[width=3.5cm,height=3.5cm]{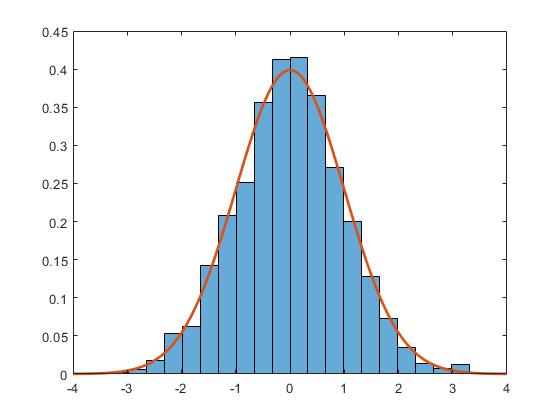}\label{w3}}
	\subfigure[Model 4]{
		\includegraphics[width=3.5cm,height=3.5cm]{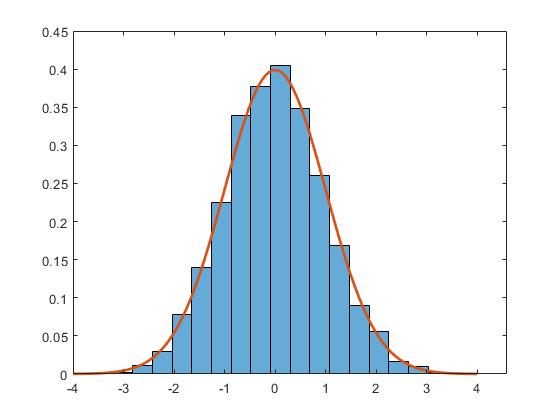}\label{w4}}
	
	\caption{
		Comparisons between empirical distributions of statistics $W$ and standard normal curves	under Models 1--4, respectively, when samples are from $Dt_1$		
	}
	\label{fig2}
\end{figure}

\begin{figure}[htbp]
	\centering
	\subfigure[Model 1]{
		\includegraphics[width=3.5cm,height=3.5cm]{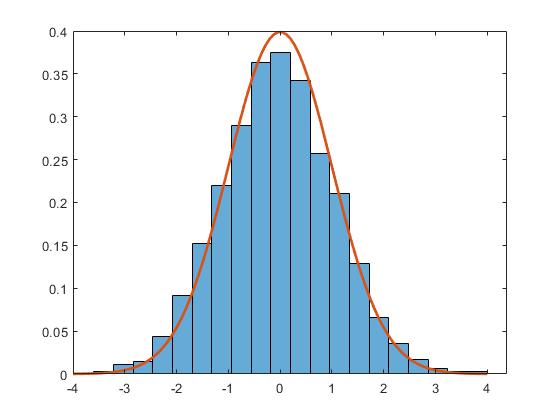}\label{v1}}
	\subfigure[Model 2]{
		\includegraphics[width=3.5cm,height=3.5cm]{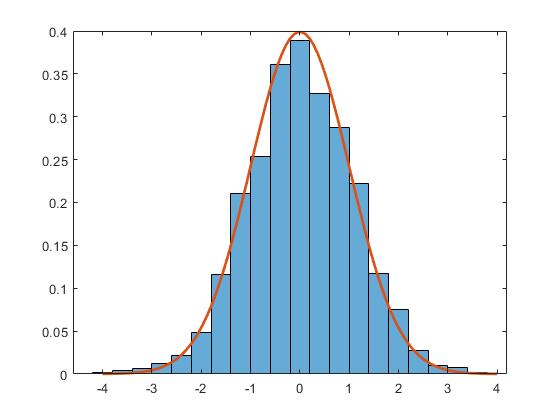}\label{v2}}
	\subfigure[Model 3]{
		\includegraphics[width=3.5cm,height=3.5cm]{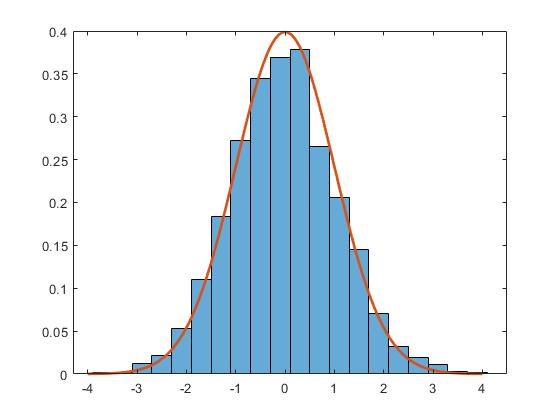}\label{v3}}
	\subfigure[Model 4]{
		\includegraphics[width=3.5cm,height=3.5cm]{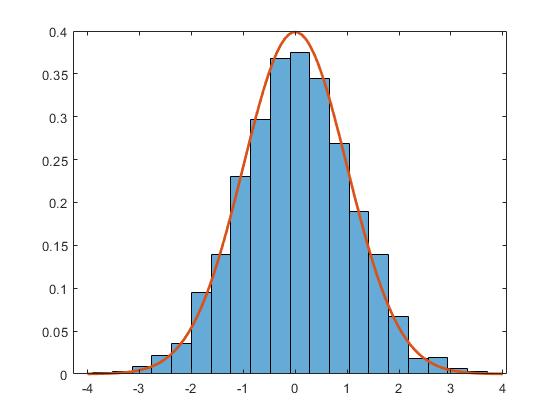}\label{v4}}
	
	\caption{
		Comparisons between empirical distributions of statistics $V$ and standard normal curves under Models 1--4, respectively, when samples are from $Dt_1$		
	}
	\label{fig3}
\end{figure}

\begin{figure}[htbp]
	\centering
	\subfigure[Model 1]{
		\includegraphics[width=3.5cm,height=3.5cm]{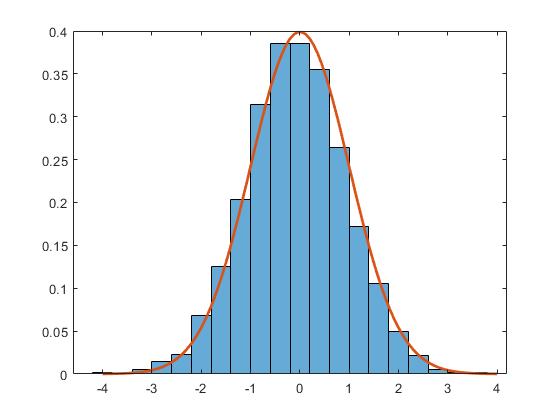}\label{u1.1}}
	\subfigure[Model 2]{
		\includegraphics[width=3.5cm,height=3.5cm]{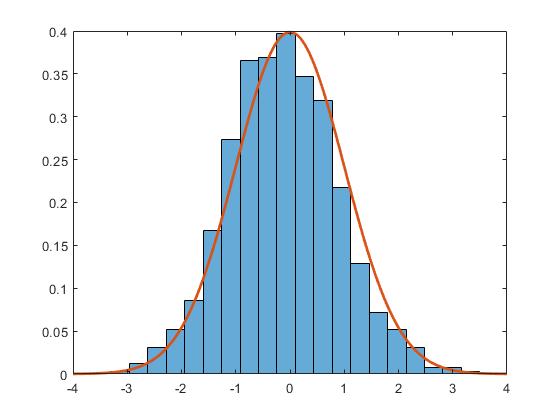}\label{u2.1}}
	\subfigure[Model 3]{
		\includegraphics[width=3.5cm,height=3.5cm]{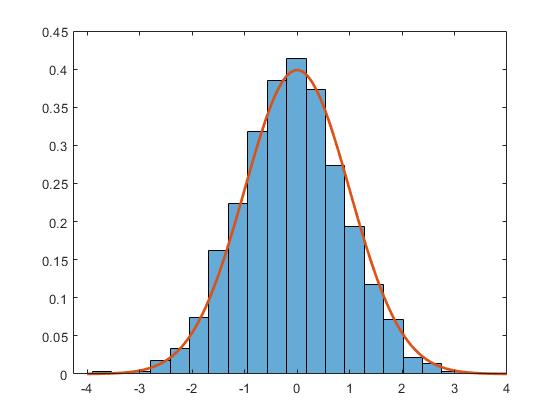}\label{u3.1}}
	\subfigure[Model 4]{ 
		\includegraphics[width=3.5cm,height=3.5cm]{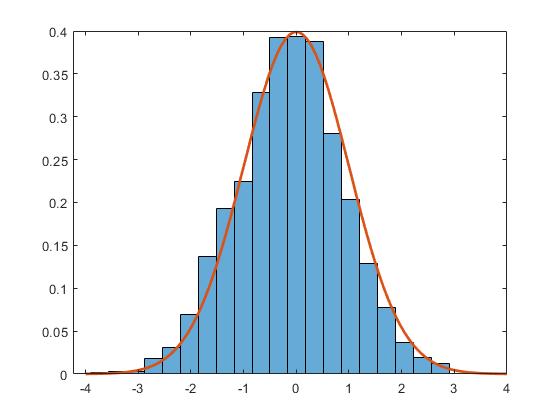}\label{u4.1}}
	
	\caption{
		Comparisons between empirical distributions of statistics $U$ and standard normal curves	under Models 1--4, respectively, when samples are from $Dt_2$	
	}
	\label{fig4}
\end{figure}

\begin{figure}[htbp]
	\centering
	\subfigure[Model 1]{
		\includegraphics[width=3.5cm,height=3.5cm]{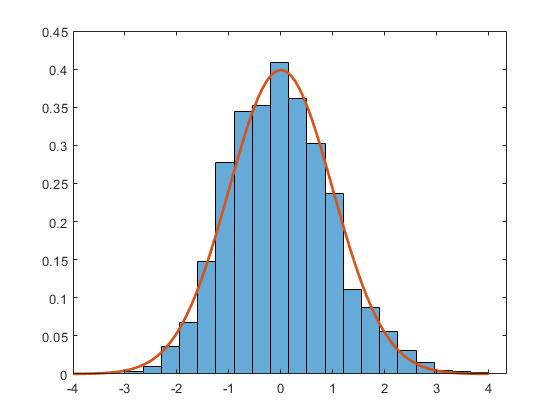}\label{w1.1}}
	\subfigure[Model 2]{
		\includegraphics[width=3.5cm,height=3.5cm]{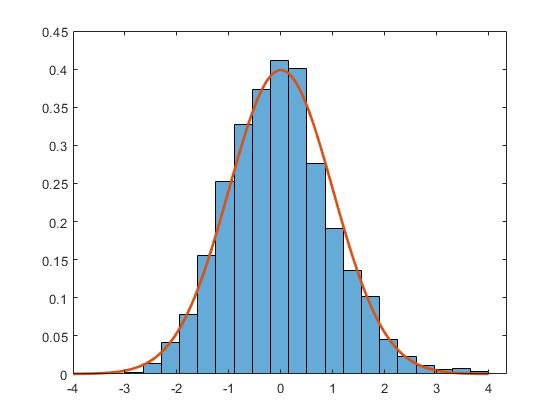}\label{w2.1}}
	\subfigure[Model 3]{
		\includegraphics[width=3.5cm,height=3.5cm]{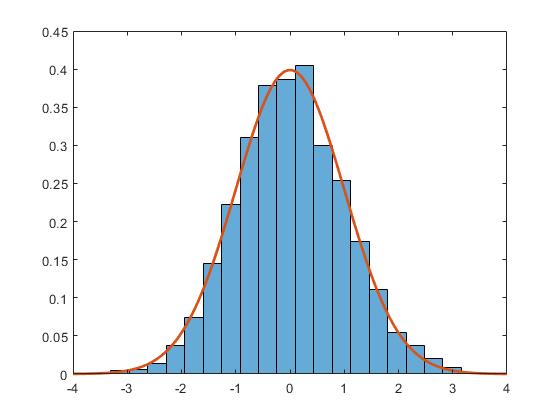}\label{w3.1}}
	\subfigure[Model 4]{
		\includegraphics[width=3.5cm,height=3.5cm]{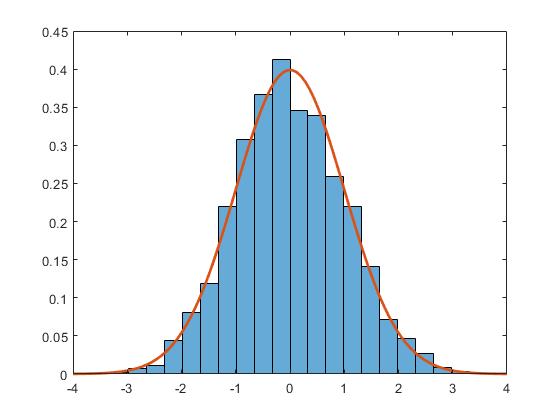}\label{w4.1}}
	
	\caption{
		Comparisons between empirical distributions of statistics $W$ and standard normal curves	under Models 1--4, respectively, when samples are from $Dt_2$		
	}
	\label{fig5}
\end{figure}

\begin{figure}[t]
	\centering
	\subfigure[Model 1]{
		\includegraphics[width=3.5cm,height=3.5cm]{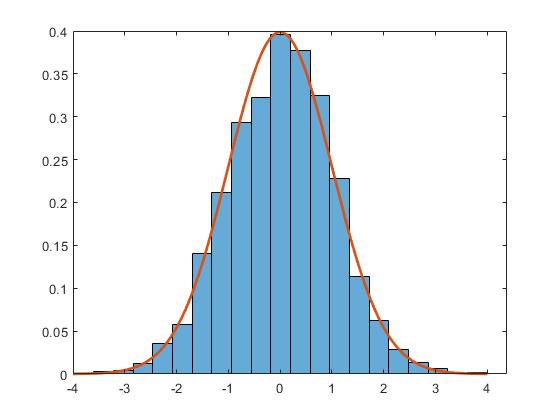}\label{v1.1}}
	\subfigure[Model 2]{
		\includegraphics[width=3.5cm,height=3.5cm]{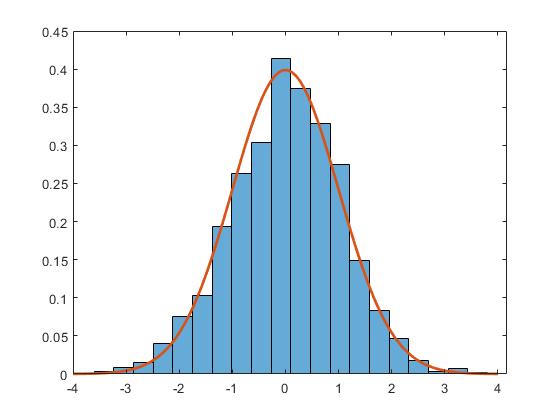}\label{v2.1}}
	\subfigure[Model 3]{
		\includegraphics[width=3.5cm,height=3.5cm]{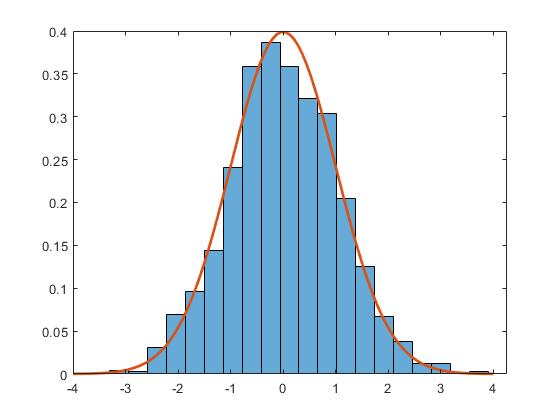}\label{v3.1}}
	\subfigure[Model 4]{
		\includegraphics[width=3.5cm,height=3.5cm]{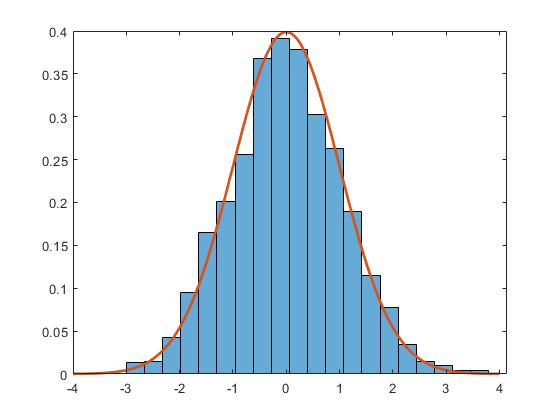}\label{v4.1}}
	
	\caption{
		Comparisons between empirical distributions of statistics $V$ and standard normal curves under Models 1--4, respectively, when samples are from $Dt_2$		
	}
	\label{fig6}
\end{figure}

Figures \ref{fig1}--\ref{fig3} show the performances of our proposed CLT (Theorems \ref{U}--\ref{V}). We compare the empirical distbributions of $U$,$W$ and $V$ with standard normal distributions (represented by red lines) under Models 1--4 when samples are from $Dt_1$. The empirical results are obtained based on 2000 replications with $p=200,n=600$. Under the same settings, when samples are from $Dt_2$, the comparisons are given in Figures \ref{fig4}--\ref{fig6}.



\end{document}